\documentclass[12pt,leqno,a4paper,notitlepage]{amsart}
\usepackage{amsmath}
\usepackage{amssymb}
\usepackage{amsfonts}
\usepackage{latexsym}
\usepackage{amscd}
\usepackage[dvips]{graphicx}

\newtheorem{teo}{Theorem}[section]
\newtheorem{prop}[teo]{Proposition}
\newtheorem{lema}[teo]{Lemma}

\newtheorem{rmk}[teo]{Remark}
\newtheorem{defn}[teo]{Definition}
\newtheorem{stp}[teo]{\newline}
\newenvironment{prova}{\textbf{Proof.}}{\hfill $\quad \Box$}

\newcommand{\Pic}{\text{\rm Pic}}

\newcommand{\TX}{{\mathcal T}(X)}
\newcommand{\STX}{{\mathcal S}{\mathcal T}(X)}
\newcommand{\TXP}{{\mathcal T}_P(X)}

\newcommand{\bigomega}{\hbox{\large $\omega$}}

\newcommand{\Lcal}{\mathcal{L}}

\newcommand{\Ical}{\mathcal{I}}

\newcommand{\Mcal}{\mathcal{M}}
\newcommand{\Ncal}{\mathcal{N}}

\newcommand{\Kcal}{\mathcal{K}}

\newcommand{\wt}{\widetilde}
\newcommand{\ol}{\overline}
\newcommand{\ox}{\otimes}
\renewcommand{\:}{\colon}

\newcommand{\IP}{{\mathbb P}}

\newcommand{\w}{\bigomega}

\newcommand{\lra}{\longrightarrow}

\renewcommand{\o}{\mathcal{O}}

\catcode`\@=11

\def\activeat#1{\csname @#1\endcsname}
\def\def@#1{\expandafter\def\csname @#1\endcsname}
{\catcode`\@=\active \gdef@{\activeat}}

\let\ssize\scriptstyle
\newdimen\ex@   \ex@.2326ex

\def\requalfill{\cleaders\hbox{$\mkern-2mu\mathord=\mkern-2mu$}\hfill
 \mkern-6mu\mathord=$}
\def\eqfill{$\m@th\mathord=\mkern-6mu\requalfill}
\def\deffill{\hbox{$:=$}$\m@th\mkern-6mu\requalfill}
\def\fiberbox{\hbox{$\vcenter{\hrule\hbox{\vrule\kern1ex
    \vbox{\kern1.2ex}\vrule}\hrule}$}}

\font\arrfont=line10
\def\Swarrow{\vcenter{\hbox{$\swarrow$\kern-.26ex
   \raise1.5ex\hbox{\arrfont\char'000}}}}

\newdimen\arrwd
\newdimen\minCDarrwd \minCDarrwd=2.5pc
       
\def\findarrwd#1#2#3{\arrwd=#3%
 \setbox\z@\hbox{$\ssize\;{#1}\;\;$}%
\setbox\@ne\hbox{$\ssize\;{#2}\;\;$}%
 \ifdim\wd\z@>\arrwd \arrwd=\wd\z@\fi
 \ifdim\wd\@ne>\arrwd \arrwd=\wd\@ne\fi}
\newdimen\arrowsp\arrowsp=0.375em
\def\findCDarrwd#1#2{\findarrwd{#1}{#2}{\minCDarrwd}
   \advance\arrwd by 2\arrowsp}
\newdimen\minarrwd 
\setbox\z@\hbox{$\longrightarrow$} \minarrwd=\wd\z@

\def\harrow#1#2#3#4{{\minarrwd=#1\minarrwd%
  \findarrwd{#2}{#3}{\minarrwd}\kern\arrowsp
   \mathrel{\mathop{\hbox to\arrwd{#4}}\limits^{#2}_{#3}}\kern\arrowsp}}
\def@]#1>#2>#3>{\harrow{#1}{#2}{#3}\rightarrowfill}
\def@>#1>#2>{\harrow1{#1}{#2}\rightarrowfill}
\def@<#1<#2<{\harrow1{#1}{#2}\leftarrowfill}
\def@={\harrow1{}{}\eqfill}
\def@:#1={\harrow1{}{}\deffill}
\def@ N#1N#2N{\vCDarrow{#1}{#2}\UpDownarrow}
\def\UpDownarrow{\uparrow\,\Big\downarrow}

\def\hookleftarrowfill{$\m@th\mathord\leftarrow\mkern-6mu%
 \cleaders\hbox{$\mkern-2mu\mathord-\mkern-2mu$}\hfill\mkern-12mu%
 \relbar\joinrel\rhook$}
\def@ l#1<#2<{\harrow1{#1}{#2}\hookleftarrowfill}

\def@'#1'#2'{\harrow1{#1}{#2}\tarrowfill}
\def\lgTo{\dimen0=\arrwd \advance\dimen0-2\arrowsp
       \hbox to\dimen0{\rightarrowfill}}
\def\smashedlgTo{\setbox0=\hbox{$\scriptstyle\lgTo$}\ht0=1.85pt
       \lower1.25pt\box0}
\def\tto{\buildrel\lgTo\over{\smashedlgTo}}
\def\tarrowfill{\hfil$\tto$\hfil}  

\def@={\ifodd\row\harrow1{}{}\eqfill
  \else\vCDarrow{}{}\Vert\fi}
\def@.{\ifodd\row\relax\harrow1{}{}\hfill
  \else\vCDarrow{}{}.\fi}
\def@|{\vCDarrow{}{}\Vert}
\def@ V#1V#2V{\vCDarrow{#1}{#2}\downarrow}
\def@ A#1A#2A{\vCDarrow{#1}{#2}\uparrow}
\def@(#1){\arrwd=\csname col\the\col\endcsname\relax
  \hbox to 0pt{\hbox to \arrwd{\hss$\vcenter{\hbox{$#1$}}$\hss}\hss}}

\def\squash#1{\setbox\z@=\hbox{$#1$}\finsm@@sh}
\def\finsm@@sh{\ifnum\row>1\ht\z@\z@\fi \dp\z@\z@ \box\z@}

\newcount\row \newcount\col \newcount\numcol \newcount\arrspan
\newdimen\vrtxhalfwd  \newbox\tempbox

\def\innernewdimen{\alloc@1\dimen\dimendef\insc@unt}
\def\measureinit{\col=1\vrtxhalfwd=0pt\arrspan=1\arrwd=0pt
  \setbox\tempbox=\hbox\bgroup$}
\def\setinit{\col=1\hbox\bgroup$\ifodd\row
  \kern\csname col1\endcsname
  \kern-\csname row\the\row col1\endcsname\fi}
\def\findvrtxhalfsum{$\egroup
 \expandafter\innernewdimen\csname row\the\row col\the\col\endcsname
 \global\csname row\the\row col\the\col\endcsname=\vrtxhalfwd
 \vrtxhalfwd=0.5\wd\tempbox
 \global\advance\csname row\the\row col\the\col\endcsname by \vrtxhalfwd
 \advance\arrwd by \csname row\the\row col\the\col\endcsname
 \divide\arrwd by \arrspan
 \loop\ifnum\col>\numcol \numcol=\col%
\expandafter\innernewdimen \csname col\the\col\endcsname
    \global\csname col\the\col\endcsname=\arrwd
  \else \ifdim\arrwd >\csname col\the\col\endcsname
     \global\csname col\the\col\endcsname=\arrwd\fi\fi
  \advance\arrspan by -1 %
  \ifnum\arrspan>0 \repeat}
\def\setCDarrow#1#2#3#4{\advance\col by 1 \arrspan=#1
   \arrwd= -\csname row\the\row col\the\col\endcsname\relax
   \loop\advance\arrwd by \csname col\the\col\endcsname
    \ifnum\arrspan>1 \advance\col by 1 \advance\arrspan by -1%
    \repeat
   \squash{\mathop{
    \hbox to\arrwd{\kern\arrowsp#4\kern\arrowsp}}\limits^{#2}_{#3}}}
\def\measureCDarrow#1#2#3#4{\findvrtxhalfsum\advance\col by 1%
  \arrspan=#1\findCDarrwd{#2}{#3}%
   \setbox\tempbox=\hbox\bgroup$}
\def\vCDarrow#1#2#3{\kern\csname col\the\col\endcsname
   \hbox to 0pt{\hss$\vcenter{\llap{$\ssize#1$}}%
    \Big#3\vcenter{\rlap{$\ssize#2$}}$\hss}\advance\col by 1}

\def\setCD{\def\harrow{\setCDarrow}%
 \def\\{$\egroup\advance\row by 1\setinit}
 \m@th\lineskip3\ex@\lineskiplimit3\ex@ \row=1\setinit}
\def\endsetCD{$\strut\egroup}
\def\measure{\bgroup
 \def\harrow{\measureCDarrow}%
 \def\\##1\\{\findvrtxhalfsum\advance\row by 2 \measureinit}%
 \row=1\numcol=0\measureinit}
\def\endmeasure{\findvrtxhalfsum\egroup}

\newbox\CDbox \newdimen\sdim

\newcount\savedcount
\def\CD#1\endCD{\savedcount=\count11%
  \measure#1\endmeasure
  \vcenter{\setCD#1\endsetCD}%
  \global\count11=\savedcount}
\catcode`\@=\active

\begin{document}

\title[{\bf Abel maps of Gorenstein curves}\hfill]
{\baselineskip=12pt \bf Abel maps of Gorenstein curves}
\author[{\bf\hfill Caporaso, Coelho, Esteves}]
{\vskip-1.7cm}

\maketitle

\

\begin{center}
{\bf Lucia Caporaso$^{\text{\bf a}}$, Juliana Coelho$^{\text{\bf b}}$ and
Eduardo Esteves$^{\text{\bf c}}$
\footnote{
Esteves was supported
by CNPq, Processos 301117/04-7 and 470761/06-7,
by CNPq/FAPERJ, Processo E-26/171.174/2003, and by the
Institut Mittag--Leffler (Djursholm, Sweden).}}
\end{center}

\vskip0.4cm

\begin{center}
$^{\text{\rm a}}$
Dipartimento di Matematica, Universit\`a Roma Tre\\
Largo San Leonardo Murialdo 1, I-00146 Roma, Italy\\
\  \\
$^{\text{\rm b}}$
Departamento de Matem\'atica, Universidade Federal Fluminense,\\
Rua M\'ario Santos Braga s/n, 24020-005 Niter\'oi, Brazil\\
\  \\
$^{\text{\rm c}}$
Rua Almeida Godinho 26, apto. 106\\
22471-140 Rio de Janeiro RJ, Brazil
\end{center}

\

\begin{abstract}
For a  Gorenstein  curve $X$ and a nonsingular point $P\in X$, we
construct Abel maps $A\:X\to J_X^1$ and  
  $A_P\:X\to J_X^0$,
where $J_X^i$ is the   moduli scheme for simple,
torsion-free, rank-1 sheaves on $X$ of degree $i$. 
The image curves  of $A$ and $A_P$ are shown to have the same
arithmetic genus of $X$. 
Also,
$A$ and $A_P$ are shown to be embeddings away from rational subcurves $L\subset X$
meeting $\ol{X-L}$ in separating nodes.
Finally 
we  establish 
a
connection with  Seshadri's 
moduli scheme $ U_X(1)$ for   semistable,
torsion-free, rank-1 
sheaves on $X$, obtaining an embedding of $A(X)$ into
$U_X(1)$.
\end{abstract}


\section{Introduction}

Fix an algebraically closed field $k$ of any characteristic
and let $X$ be a connected,   projective
curve over $k$.
If $X$ is   smooth, there is, for each
integer $d\geq 1$, a natural map,
$$
A^d\: X^d \longrightarrow\Pic^dX,
$$
with $\Pic^dX$ denoting
the Picard scheme parameterizing line bundles of
degree $d$ on $X$; the map sends $(P_1,\dots,P_d)$ to
$[\o_X(P_1+\cdots+P_d)]$. 
Of course
$A^d$ factors through a map
$$
A^{(d)}\: X^{(d)}\longrightarrow\Pic^dX,
$$
where $X^{(d)}$ is the $d$th symmetric product of $X$.
The map $A^{(d)}$ is called the degree-$d$ \emph{Abel map} of $X$.

Classical variants are the degree-$d$
Abel maps with base point $P\in X$,
$$
A^d_P\: X^{(d)}\longrightarrow\Pic^0X,
$$
which is simply $A^d$ composed with the translation map
$\Pic^dX\to\Pic^0X$, taking 
$[\Lcal]$ to $[\Lcal\ox\o_X(-dP)]$,
and the compositions of $A^{(d)}$ and $A^{(d)}_P$ with
the duality isomorphisms $\lambda^d\:\Pic^d_X\to\Pic^{-d}_X$, taking
$[\Lcal]$ to $[\Lcal^*]$.

Much of the geometry of $X$ is encoded in the Abel maps,
since their fibers are the 
complete linear systems of $X$. For instance, the gonality of $X$ is
the smallest integer $d$ such that $A^{(d)}$ is not an embedding. In
particular, $X$ is hyperelliptic if and only if $A^{(2)}$ is not an
embedding.

The Abel maps behave naturally in families of smooth curves.
As smooth curves degenerate to singular ones, we would like to
understand how the Abel maps degenerate. So,
how to construct Abel maps for singular curves in a natural way?

If $X$ is 
integral,
Altman and Kleiman
\cite{AK} defined, for each 
$d\geq 1$, a natural map
$$
\beta^{(d)}\:\text{Hilb}^dX\longrightarrow J^{-d}_X,
$$
where $\text{Hilb}^dX$ is the Hibert scheme of $X$ parameterizing
length-$d$ subschemes, and $J^{-d}_X$ is the
compactified Jacobian parameterizing torsion-free, rank-1 sheaves
of degree $-d$ on $X$; the map sends $[Y]$ to
$[\Ical_{Y/X}]$. Again, the fibers of $\beta^{(d)}$ are projective
spaces. And, if $X$ is smooth, then
$\text{Hilb}^dX=X^{(d)}$, $J^{-d}_X=\Pic^{-d}X$ and
$\beta^{(d)}=\lambda^d\circ A^{(d)}$.

On the other hand, if a curve is reducible, the situation is more
complex. The 
current knowledge
is concentrated 
on 
the two extremes:
$d=1$ and $d=g-1$.
For $d=g-1$, the image of $A^d$ 
turns out to be
the 
theta divisor. 
For work  extending the construction 
and the properties
of the theta divisor
to singular curves
we refer the reader to
\cite{So} and \cite{E97} for 
irreducible curves, and
to
the more recent \cite{A04} and \cite{captheta1} 
for 
nodal, possibly reducible, curves.

As for $d=1$,  Edixhoven \cite{Edix98} constructed and studied
rational Abel maps of nodal curves to N\'eron
models. As N\'eron models are
seldom complete, 
his maps 
are
not defined everywhere. In
\cite{capoed} the compactifications of N\'eron models 
and Picard schemes,  constructed in \cite{cner}, 
are used;   it is shown
that, if $X$ is  stable, 
there exists a
globally defined map $\ol{\alpha^1_X}\:X\to\ol{P^1_X}$,
where $\ol{P^1_X}$ 
is the compactified Picard scheme parameterizing
equivalence
classes of degree-1 ``semibalanced'' line bundles on semistable
curves having $X$ as a stable model. 

In this paper we extend the construction of $\ol{\alpha^1_X}$ to any
G-stable curve $X$, that is, to any 
reduced
curve $X$ 
with Gorenstein singularities whose dualizing
sheaf is ample. 
Also, we
describe 
the image and the fibers of our
Abel map.
More precisely, for any 
reduced 
curve $X$ 
with Gorenstein singularities, we consider
the fine moduli schemes $J^d_X$, parameterizing simple,
torsion-free, rank-1 sheaves of degree $d$ on $X$, for all
integers $d$.
And we construct a map (cf.~\ref{twembpre}):
$$
A\: X\lra J_X^1.
$$
The schemes $J^d_X$ are quite large, not even Noetherian, but have
an open subscheme of finite type, $J^{d,ss}_X$, parameterizing
semistable sheaves (cf. \ref{finejac}). The map
$A$ is constructed in such a way that $A(X)\subseteq J^{1,ss}_X$,
if $X$ is G-stable (cf. Theorem~\ref{twemb}).

If $X$ has no separating nodes
(nodes whose removal disconnects the curve), 
then $A$ sends $Q$ to
$[\Ical_{Q/X}^*]$. In particular, if $X$ is smooth we
recover the classical degree-1 Abel map $A^{(1)}$. 
On the other hand,
if $X$ 
does admit a
separating node $N$, then $\Ical_{N/X}^*$ is not simple,
and thus not parameterized by $J_X^1$. So, for each $Q\in X$
we create a new sheaf $\Ical_Q^1$ out of $\Ical_{Q/X}^*$,
by tensoring the latter with 
suitable so-called ``twisters",
along the same lines of what was done in \cite{capoed}, and let
$A$ send $Q$ to $[\Ical_Q^1]$.

In Theorem \ref{cptilde} we prove that $A$ contracts every
smooth rational subcurve  $L\subseteq X$ 
meeting its complementary curve in separating nodes of $X$,
and  $A$ is an embedding off these subcurves. 
Also, $A(X)$ has the same arithmetic
genus of $X$ and its singularities are those of $X$, 
together with ordinary
singularities with linearly independent tangent lines.

Unfortunately, 
the schemes $J^{d,ss}_X$ 
are not, in
general, separated. To get a separated scheme, two alternatives are
possible:  
to
use either smaller schemes or  quotient schemes.
We consider both.

For each 
nonsingular point $P\in X$
the scheme $J^{d,ss}_X$ has an open subscheme
$J^{d,P}_X$, parameterizing sheaves that are
$P$-quasistable (cf. \ref{finejac}),
which is 
projective over $k$. If $P$ is suitably
chosen, and $X$ is G-stable, then $A(X)\subseteq J^{1,P}_X$ by
our Theorem \ref{twemb}.

On the other hand, we consider Seshadri's coarse moduli schemes $U_X(d)$ 
for equivalence classes of semistable, torsion-free, rank-1 sheaves on
$X$ of degree $d$ (cf.~\ref{seshjac}). There are natural maps
$\Phi^d\:J^{d,ss}_X\to U_X(d)$, taking a semistable sheaf to its
class. Our Theorem \ref{seshsep} says that $\Phi^1$ restricts
to an embedding on $A(X)$. 

We 
also
construct Abel maps with base points. More precisely,
for each nonsingular  point  $P\in X$
we construct a
natural map (cf. \ref{twptembpre}),
$$
A_P\:X\lra J_X^0,
$$
with image in $J_X^{0,P}$ (cf. Theorem \ref{twptemb}).
If $X$ has no separating nodes,
then $A_P$ sends $Q$ to $[\Ical_{Q/X}\ox\o_X(P)]$. So, if $X$ is
smooth then $A_P=\lambda^0\circ A_P^{(1)}$.
If $X$ has separating nodes, $A_P$ is constructed
with the help of twisters, as done for $A$.

The map $A_P$ has the same description as $A$. In fact, for a
suitably chosen $P$, we may view $A$ as the composition of
$A_P$ with the duality map $J_X^0\to J_X^0$, taking $[\Ical]$ to
$[\Ical^*]$ and the translation map $J_X^0\to J_X^1$, taking
$[\Ical]$ to $[\Ical\ox\o_X(P)]$. So, up to these isomorphisms,
$A$ may be viewed as one of the $A_P$. 
(In fact, everywhere in the paper we prove properties first for $A_P$ and then
extend them to $A$.)

The biggest difference between $A$ and the $A_P$ comes 
when we consider their composition
to Seshadri's moduli schemes: $\Phi^0\circ A_P$ may actually
collapse components of $X$ that were not collapsed by $A_P$, as
we point out in Remark \ref{0bad}.

We conclude with
a few   comments 
about closely related questions and further  developements. 
First, $A$ is not always
a natural map. In fact, if $X$ has a ``splitting" node  
(i.e. a
separating node splitting the curve in two equal genus subcurves)
then $A$ depends on the choice
of one of these subcurves (cf.~\ref{twembpre}). The lack of
naturality is a major hurdle
to extend our construction to families of curves.
In this respect,
the 
map to Seshadri's  moduli space, 
$\Phi^1\circ A$, 
looks
more natural, as 
it is
independent of the above choice  by
Theorem \ref{seshsep}.

On the other hand, 
the maps $A_P$ are 
natural, and it
seems possible to extend their construction to families of pointed
curves, whereas the compositions $\Phi^0\circ A_P$ do not behave
well. 
We hope to deal with Abel maps for families  in the future.

Second, we don't 
treat
higher degree Abel maps. 
It seems
possible  to define them  not on $X^d$, $X^{(d)}$ or 
$\text{Hilb}^d_X$, but on blowups of 
them. Very little is known, apart from  the case of
 the degree-2 Abel
map for a nodal curve with two components meeting at two points,
constructed in  \cite{Co06}.

Here is a layout of the paper: In Section 2, we
introduce the moduli schemes $J_X^{d,ss}$ and $U_X(d)$,
and the quotient maps $\Phi^d\:J_X^{d,ss}\to U_X(d)$.
In Section~3, we construct the Abel
maps $A$ and $A_P$ when $X$ has no separating nodes.
In Section~4, we construct the maps $A_P$ in general, and in
Section 5 we do the same for the map $A$. In Section 6, we prove
properties of the maps $A$ and $A_P$, describing their images and
fibers. Finally, in Section 7 we show that $\Phi^1$ restricts
to an embedding on $A(X)$.

\section{Compactified Jacobians}

All schemes are assumed to be locally of finite type over 
a fixed algebraically closed field  
$k$.
A point of a scheme means a closed point.
A \emph{curve} is a  reduced, projective scheme of
pure dimension 1 over $k$. If $Y$ is a curve, we let
$g_Y:=1-\chi(\o_Y)$ and call $g_Y$ the (arithmetic) \emph{genus} of
$Y$.

\emph{Throughout the paper,
$X$ denotes a connected curve, $\w$ its dualizing sheaf,
$g$ its (arithmetic) genus and $P$ a
point on the nonsingular locus of $X$.}

\begin{stp}\label{setup}\setcounter{equation}{0}\rm
(\emph{Preliminaries})
A reduced union of irreducible components of $X$, connected or
not, is called a \emph{subcurve}.
If $Y$ is a proper subcurve of $X$, let $Y'$ denote the
\emph{complementary subcurve}, that is, the reduced union
of all the irreducible components of $X$ not contained in $Y$.
The intersection $Y\cap Y'$ is a finite scheme; let
$\delta_Y$ denote its length.
Since $X$ is connected, $\delta_Y>0$.
Also, observe that
\begin{equation}
\label{subgenus}
g_Y\leq g.
\end{equation}

Let $\Ical$ be a coherent sheaf on $X$. We say that $\Ical$ is
\emph{torsion-free} if its associated points
are generic points of $X$. We say that $\Ical$ is \emph{of rank $1$}
if $\Ical$ is invertible on a dense open subset of $X$. And we say
that $\Ical$ is \emph{simple} if $\text{End}(\Ical)=k$. Each
line bundle on $X$ is torsion-free of rank 1 and simple.

Suppose $\Ical$ is torsion-free of rank 1. We call
$\deg(\Ical):=\chi(\Ical)-\chi(\o_X)$ the \emph{degree} of $\Ical$.
For each vector bundle $F$ on $X$,
\begin{equation}\label{IoxF}
\chi(\Ical\ox F)=\text{rk}(F)\chi(\Ical)+\deg(F)=
\text{rk}(F)(\deg(\Ical)+1-g)+\deg(F).
\end{equation}

For each subcurve $Y$ of $X$, let $\Ical_Y$ denote the restriction of
$\Ical$ to $Y$ modulo torsion, that is, the image of the natural map
$$
\Ical|_Y\lra\bigoplus_{i=1}^m  (\Ical|_Y)_{\xi_i},
$$
where $\xi_1,\dots,\xi_m$ are the generic points of $Y$. We let
$\deg_Y(\Ical)$ denote the degree of $\Ical_Y$, that is,
$\deg_Y(\Ical):=\chi(\Ical_Y)-\chi(\o_Y)$.

Let $Y$ be a proper subcurve of $X$.
By the defining property of the dualizing sheaf,
the kernel of the restriction
map $\w\to\w|_{Y'}$ is the dualizing sheaf $\w_Y$ of $Y$.
Suppose $X$ is Gorenstein. Then
$\w$ is a line bundle of degree $2g-2$,
and it follows that
$$
\chi(\w|_Y)=\chi(\w_Y)+\delta_Y.
$$
Thus, by duality,
\begin{equation}\label{dwY1}
\deg(\w|_Y)=\chi(\w_Y)+\delta_Y-\chi(\o_Y)=
-2\chi(\o_Y)+\delta_Y=2g_Y-2+\delta_Y.
\end{equation}
\end{stp}

\begin{defn}\label{Gstable}\setcounter{equation}{0}\rm
The curve $X$ is called
\emph{G-stable} if $X$ is Gorenstein of genus $g\geq 2$, and
does not contain any smooth rational
component $L$ with $\delta_L\geq 2$. Equivalently, using
\eqref{dwY1}, the curve $X$ is G-stable if $X$ is Gorenstein and
$\w$ is ample.
\end{defn}

If the singularities of $X$ are (ordinary) nodes, then $X$ is G-stable
if and only if it is stable, in the sense of Deligne and Mumford.

\begin{stp}\label{cpol}\setcounter{equation}{0}\rm
(\emph{Semistable sheaves})
Let $F$ be a vector bundle and
$\Ical$ a torsion-free, rank-1 sheaf of degree
$d$ on the curve $X$.
We call $\Ical$ \emph{semistable} with respect to
$F$ if
\begin{enumerate}
\item\label{ss1}
$\chi(\Ical\ox F)=0$ and
\item\label{ss2}
$\chi(\Ncal\ox F)\geq 0$ for each nonzero quotient $\Ncal$ of
$\Ical$ different from $\Ical$.
\end{enumerate}
We call $\Ical$ \emph{stable} with respect
to $F$ if the inequality in Condition \ref{ss2} is
always strict.

By \eqref{IoxF}, Condition \ref{ss1} is verified
if the \emph{slope}
$\mu(F):=\deg(F)/\text{rk}(F)$ satisfies $\mu(F)=g-1-d$.
As for Condition \ref{ss2}, observe that all nonzero torsion-free
quotients of $\Ical$ are of the form $\Ical_Y$ for a subcurve
$Y\subseteq X$. So, Condition~\ref{ss2} holds
if and only if
\begin{equation}\label{ssrk1}
\chi(\Ical_Y\ox F|_Y)\geq 0
\end{equation}
for each proper subcurve $Y$ of $X$. For stability,
we require strict inequalities.

We say that $\Ical$ is \emph{$P$-quasistable} with respect to
$F$ if Inequality \eqref{ssrk1} holds for each proper subcurve
$Y\subset X$, with equality only if $P\not\in Y$.
Clearly, this notion depends only on
which component of $X$ the point $P$ lies.

Suppose $X$ is Gorenstein. Define a vector bundle $E_d$ on $X$
as follows: If $g\geq 2$,
\begin{equation}
\label{Ed}
E_d:=\o_X^{\oplus 2g-3}\oplus\w^{\ox g-1-d};
\end{equation}
if $g=0$, set $E_d:=\o_X\oplus\w^{\ox d+1}$; and if $g=1$, define
$E_d$ only if $d=0$, setting $E_0:=\o_X$. We call $E_d$ the
\emph{canonical $d$-polarization} of $X$.
Notice that $E_d$ has slope $\mu(E_d)=g-1-d$.

We say that $\Ical$ is (canonically)
stable, semistable or $P$-quasistable if
$\Ical$ is so with respect to $E_d$. If $g\geq 2$ then, for each
subcurve $Y\subseteq X$,
\begin{align*}
\chi(\Ical_Y\ox E_d|_Y)=&(2g-2)\chi(\Ical_Y)+\deg(E_d|_Y)\\
=&(2g-2)\deg_Y(\Ical)+(2g-2)\chi(\o_Y)+(g-1-d)\deg(\w|_Y)\\
=&(2g-2)\deg_Y(\Ical)+(1-g)(\deg(\w|_Y)-\delta_Y)+(g-1-d)\deg(\w|_Y)\\
=&(2g-2)\deg_Y(\Ical)+(g-1)\delta_Y-d\deg(\w|_Y),
\end{align*}
where we used \eqref{IoxF} and \eqref{dwY1}.
Thus, $\Ical$ is semistable
if and only if
\begin{equation}\label{ss}
\deg_Y(\Ical)\geq d\Big(\frac{\deg_Y(\w)}{2g-2}\Big)
-\frac{\delta_Y}{2}
\end{equation}
for each proper subcurve $Y$ of $X$.
If $g=0$, then an analogous computation can be done, and the
condition for semistability is the same. Finally,
if $g=1$ and $d=0$,
then $\Ical$ is semistable if and only if
\begin{equation}\label{sss}
\deg_Y(\Ical)\geq-\frac{\delta_Y}{2}
\end{equation}
for each proper subcurve $Y$ of $X$. Notice that
\eqref{ss} and \eqref{sss} are equal conditions if $d=0$.
We leave it to the reader to formulate the analogous conditions
for when $\Ical$ is stable or $P$-quasistable.
\end{stp}

\begin{stp}\label{finejac}\setcounter{equation}{0}\rm
(\emph{The fine compactified Jacobians}) There exists a scheme $J_X$
parameterizing torsion-free, rank-1,
simple sheaves on the curve $X$; see \cite{ed01} Thm. B, p. 3048.
More precisely, given a scheme $T$, a $T$-flat coherent sheaf $\Ical$ on
$X\times T$ is called torsion-free (resp. rank-1, resp. simple) on
$X\times T/T$ if $\Ical|_{X\times t}$ is torsion-free (resp. rank-1, resp.
simple) for every $t\in T$. The scheme $J_X$ represents
the functor that associates to each scheme $T$ the set of
torsion-free, rank-1 simple sheaves on $X\times T/T$
modulo equivalence $\sim$. Two such sheaves $\Ical_1$ and $\Ical_2$
are called equivalent, $\Ical_1\sim\Ical_2$,
if there is a line bundle $N$ on $T$ such that
$\Ical_1\cong\Ical_2\ox p_2^*N$, where $p_2\:X\times T\to T$ is the
projection map.

If $T$ is a connected scheme, and $\Ical$ is a torsion-free,
rank-1 sheaf on $X\times T/T$, then $d:=\deg(\Ical|_{X\times t})$
does not depend on the choice of $t\in T$; we say that
$\Ical$ is a degree-$d$ sheaf on $X\times T/T$. Then
there is a natural decomposition
$$
J_X=\coprod_{d\in\text{\bf Z}}J^d_X,
$$
where $J^d_X$ is the subscheme of $J_X$
parameterizing degree-$d$ sheaves.
The schemes $J^d_X$ are universally closed
over $k$; see \cite{ed01} Thm. 32, (2), p.~3068. However,
in general, the $J^d_X$ are neither of finite type nor separated over
$k$.

Let $F$ be a vector bundle on $X$ with integer slope, and set
$d:=g-1-\mu(F)$. By \cite{ed01} Prop. 34, p. 3071, the subschemes
$J^{ss}_F$ (resp. $J^s_F$, resp. $J^P_F$) of $J^d_X$
parameterizing simple and semistable (resp. stable, resp.
$P$-quasistable)
sheaves on $X$ with respect to $F$ are open.
If $X$ is Gorenstein and $F=E_d$
(the canonical $d$-polarization defined in \ref{cpol}),
we write
$$
J^{d,ss}_X:=J^{ss}_{E_d},\quad J^{d,s}_X:=J^s_{E_d}\quad\text{and}\quad
J^{d,P}_X:=J^P_{E_d}.
$$

By \cite{ed01} Thm. A, p. 3047,
$J^{ss}_F$ is of finite type and
universally closed, $J^s_F$ is separated and
$J^P_F$ is complete over $k$. Actually, $J^P_F$ is projective; see
\cite{CoE} Prop.~2.4.
\end{stp}

\begin{stp}\label{Sequiv}\setcounter{equation}{0}\rm
(\emph{The S-equivalence}) Let $F$ be a vector bundle on $X$. For
each semistable sheaf $\Ical$ on $X$ with respect to $F$, there is a
maximal filtration
$$
\emptyset\subsetneqq Y_1\subsetneqq Y_2\subsetneqq\cdots\subsetneqq
Y_{q-1}\subsetneqq X
$$
of $X$ by subcurves $Y_i$ such that $\chi(\Ical_{Y_i}\ox F|_{Y_i})=0$
for each $i=1,\dots,q$, which we call a
\emph{Jordan--H\"older filtration}. There may be many
Jordan--H\"older filtrations associated to $\Ical$, but
the collection of subcurves
$$
\mathfrak S(\mathcal I):=\{Y_1,\ol{Y_2-Y_1},\dots,\ol{Y_{q-1}-Y_{q-2}},
\ol{X-Y_q}\}
$$
and the isomorphism class of the sheaf
$$
\text{Gr}(\Ical):=\Ical_{Y_1}\oplus\text{Ker}(\Ical_{Y_2}\to\Ical_{Y_1})
\oplus\cdots\oplus\text{Ker}(\Ical_{Y_{q-1}}\to\Ical_{Y_{q-2}})\oplus
\text{Ker}(\Ical\to\Ical_{Y_q})
$$
depend only on $\Ical$, by the Jordan--H\"older Theorem.

We say that two semistable
sheaves $\Ical$ and $\Kcal$ on $X$ are \emph{S-equivalent} if
$\mathfrak S(\Ical)=\mathfrak S(\Kcal)$ and
$\text{Gr}(\Ical)\cong\text{Gr}(\Kcal)$.

(For a higher rank semistable sheaf, a Jordan--H\"older filtration
is a filtration of the sheaf. However, in rank 1, this filtration
is induced by a filtration of $X$ as above.)
\end{stp}

\begin{stp}\label{seshjac}\setcounter{equation}{0}\rm
(\emph{The coarse compactified Jacobians})
Let $X_1,\dots,X_n$ be
the irreducible components of the curve $X$,
and $\mathfrak a=(a_1,\dots,a_n)$ a $n$-tuple of
positive rational numbers summing up to 1. For each
subcurve $Y\subseteq X$, set $ a_Y:=\sum_{X_i\subseteq Y}a_i$.
According to Seshadri \cite{sesh} D\'ef. 9 and
Remarques on p. 153, a
torsion-free, rank-1 sheaf
$\Ical$ on $X$ is \emph{$\mathfrak a$-semistable} if
$$
\chi(\Ical_Y)\geq a_Y\chi(I)
$$
for each proper subcurve $Y$ of $X$. Also, $\Ical$ is
\emph{$\mathfrak a$-stable} if the inequalities are strict.

Seshadri's notions of semistability and stability are encompassed by 
ours. More
precisely, for each integer $d$ there is a vector bundle $F$ on $X$
such that $\mathfrak a$-semistability (resp. $\mathfrak a$-stability)
for degree-$d$, torsion-free, rank-1 sheaves is equivalent to
semistability (resp. stability) with respect to $F$; see
\cite{edsep} Obs. 13, p. 584. In fact,
any vector bundle $F$ on $X$
such that
\begin{equation}\label{Ea}
\mu(F|_{X_i})=a_i(g-1-d)\quad\text{for each $i=1,\dots n$}
\end{equation}
has this property.

Two $\mathfrak a$-semistable sheaves are called S-equivalent if they
are S-equivalent in the sense of \ref{Sequiv} for a (and hence any)
vector bundle $F$ on $X$ satisfying \eqref{Ea}.

In \cite{sesh} Thm. 15, p. 155, Seshadri shows that there is a
scheme $U_X(\mathfrak a,d)$ corepresenting the functor
that associates to each scheme $T$ the set of
$T$-flat coherent sheaves $\Ical$ on $X\times T$ such that
$\Ical|_{X\times t}$ is $\mathfrak a$-semistable and of degree $d$
for every  $t\in T$, modulo the same equivalence $\sim$
of \ref{finejac}. Furthermore, $U_X(\mathfrak a,d)$
is projective and parameterizes S-equivalence classes of
$\mathfrak a$-semistable sheaves.

Let $F$ be any vector bundle on $X$ satisfying \eqref{Ea}, and set
$J_X(\mathfrak a,d):=J^{ss}_F$.
(The particular $F$ is irrelevant.) Since $J^{ss}_F$
is a fine moduli space, there is a naturally induced morphism
\begin{equation}\label{Smap}
\Phi_{\mathfrak a}^d\: J_X(\mathfrak a,d)\longrightarrow U_X(\mathfrak a,d)
\end{equation}
sending $[\Ical]$ to the S-equivalence class of $\Ical$.
We call $\Phi_{\mathfrak a}^d$ the \emph{S-map}.

If $X$ is G-stable (cf. \ref{Gstable}),
let $U_X(d):=U_X(\mathfrak a,d)$, where
$$
a_i:=\frac{\deg_{X_i}(\w)}{2g-2}\quad\text{for each $i=1,\dots,n$}.
$$
(Notice that $a_1+\cdots+a_n=1$ because $X$ is Gorenstein,
and the $a_i$ are positive because $\w$ is ample.)
Since, for each integer $i$,
$$
\mu(E_d|_{X_i})=a_i(g-1-d),
$$
where $E_d$ is the canonical $d$-polarization of $X$
(cf. \ref{cpol}), the S-map
\eqref{Smap} becomes
\begin{equation}\label{Sdmap}
\Phi^d\:J^{d,ss}_X\lra U_X(d).
\end{equation}
\end{stp}

\section{Abel maps}

\emph{Assume from now on until the end of the paper that
$X$ is Gorenstein.}

\begin{defn}\setcounter{equation}{0}\rm
A \emph{separating node} of the curve $X$ is a
point $N$ for which there is a subcurve $Z$
such that $\delta_Z=1$ and $Z\cap Z'=\{N\}$.
\end{defn}

 Being $X$  Gorenstein, a separating node is indeed a node,
 by \cite{cat}, Prop. 1.10, p. 59.

\begin{stp}\label{3.1}\setcounter{equation}{0}\rm
(\emph{Degree-$0$ Abel maps})
For each point $Q$ on the curve $X$, its
sheaf of ideals $\mathfrak m_Q$ is torsion-free of rank 1 and
degree $-1$. Also, if $Q$ is not a separating node,
$\mathfrak m_Q$ is simple, as it follows from the discussion in
\cite{ed01}, Ex. 38, p. 3073.

Let $\Ical_\Delta$ be the ideal sheaf of the diagonal
$\Delta\subset X\times X$, and put
$$
\Ical:=\Ical_\Delta\ox p_1^*\o_X(P),
$$
where $p_1\:X\times X\to X$ is the first projection.
The sheaf $\Ical$ is flat over $X$ and, for each $Q\in X$,
$$
\Ical|_{X\times Q}=\mathfrak m_Q\ox\o_X(P).
$$
If $X$ is free from separating nodes, then
$\Ical$ defines a morphism
\begin{equation}
\label{AP}
A_P\: X\lra J^0_X;\  \  \  \  Q\mapsto [\mathfrak m_Q\ox\o_X(P)].
\end{equation}
We call $A_P$ the \emph{degree-$0$ Abel map} of $X$
with base $P$.
\end{stp}

\begin{prop}
\setcounter{equation}{0}\label{pointabel}
Assume that the curve $X$ is free from separating nodes. Then:
\begin{enumerate}
\item $A_P(X)\subseteq J^{0,P}_X$.
\item If $X\not\cong\IP^1$ then $A_P$ is an embedding.
\end{enumerate}
\end{prop}

\begin{prova} For each $Q\in X$, its sheaf of ideals $\mathfrak m_Q$ 
satisfies
$$
\deg_Y(\mathfrak m_Q)=\begin{cases}
-1,&\text{if }Q\in Y,\\
0,&\text{if }Q\not\in Y,
\end{cases}
$$
for each subcurve $Y\subseteq X$. If $Y$ is proper, then
$\delta_Y>1$ by hypothesis, and hence
$\deg_Y(\mathfrak m_Q\ox\o_X(P))\geq-\delta_Y/2$, with
equality only if $P\not\in Y$. So
$\mathfrak m_Q\ox\o_X(P)$ is $P$-quasistable, showing that
$A_P(X)\subseteq J_X^{0,P}$.

Assume $X\not\cong\IP^1$.
It remains to show that $A_P$ is an embedding.
Since $X$ is complete and $J^{0,P}_X$ is separated,
the induced map $X\to J^{0,P}_X$ is proper.
Thus, we need only show that $A_P$
separates points and tangent vectors. Equivalently,
we need only show that every fiber of
$A_P$ is either empty or schematically a
point.

Let $Q\in X$ and put
$L:=A_P^{-1}([\mathfrak m_Q\ox\o_X(P)])$. From
\cite{AK} Lemma 5.17, p.~88, it follows that $L$ is
isomorphic to an open subscheme of the projective space
$$
{\mathbb P}(\text{Hom}_X(\mathfrak m_Q,\o_X)),
$$
the open subscheme parameterizing injective homomorphisms.
However, since $A_P$ is proper over $J^{0,P}_X$, the fiber
$L$ is complete, and
thus $L$ is a projective space.

We need to show that $L$ is a
point. Suppose otherwise, by contradiction. Thus,
since $L$ has dimension at most 1, we have $L\cong\IP^1$.

Let $Q_1$ and $Q_2$
be distinct points of $L$ on the nonsingular locus of
$X$. Since
$L$ is a fiber of $A_P$, we have an isomorphism
$\mathfrak m_{Q_1}\to\mathfrak m_{Q_2}$. This isomorphism
is given by multiplication by a rational function $h$ of $X$,
whose only pole is $Q_1$ and whose only zero is $Q_2$, both
with order 1. The function $h$ is constant on all components
of $X$ other than $L$, because $h$ has no zeros or poles there.
Let $Z:=L'$. Since $X$ is not isomorphic
to $\IP^1$, we have $Z\neq\emptyset$. We claim that
$L$ intersects $Z$ transversally. In fact, if
$L$ intersected $Z$ nontransversally at a point
$R$, then $h|_L-h(R)$ would vanish at $R$
with order at least $2$.
This is not possible
because $h|_L$ has degree 1.

Let $Z_1,\dots,Z_q$ denote the connected components of $Z$.
Since $X$ is connected, each $Z_i$
intersects $L$. If $\#(Z_i\cap L)=1$,
then we would have $\delta_{Z_i}=1$, as we
already know that $Z$ intersects $L$ transversally.
Thus, each $Z_i$ intersects $L$ in at least two points.
But then $h|_L$ takes the same value on these two points of
$L$. This is again not possible because
$h|_L$ has degree 1. We have reached a contradiction.
\end{prova}

\begin{stp}\setcounter{equation}{0}\rm
(\emph{Degree-$1$ Abel maps})
As in \ref{3.1},
let $\Ical_\Delta$ be the ideal sheaf of
the diagonal $\Delta\subset X\times X$.
Then $\Ical_\Delta$ is flat over $X$ and, for each $Q\in X$, the
restriction $\Ical_\Delta|_{X\times Q}$ is isomorphic to the sheaf of
ideals $\mathfrak m_Q$. Since $X$ is Gorenstein,
the dual sheaf
$$
\Ical^*_\Delta:=Hom_{X\times X}(\Ical_\Delta,\o_{X\times X})
$$
is also flat over $X$, and
$\Ical^*_\Delta|_{X\times Q}\cong\mathfrak m^*_Q$
for each $Q\in X$.

As mentioned in \ref{3.1}, the sheaf
$\mathfrak m_Q$ is simple if $X$ is free from
separating nodes. Since $X$ is Gorenstein,
$$
\text{Hom}_X(\mathfrak m_Q^*,\mathfrak m_Q^*)=
\text{Hom}_X(\mathfrak m_Q,\mathfrak m_Q).
$$
So, if $X$ is free from separating nodes,
$\Ical^*_\Delta|_{X\times Q}$
is simple for every $Q\in X$, and thus
$\Ical^*_\Delta$ defines a morphism
\begin{equation}\label{Anotwist}
A\: X\lra J^1_X; \  \  \  Q\mapsto[\mathfrak m_Q^*].
\end{equation}
We call $A$ the \emph{degree-$1$ Abel map} of $X$.
\end{stp}

\begin{prop}
\label{abel}\setcounter{equation}{0}
Assume that the curve $X$ is free from separating nodes. Then:
\begin{enumerate}
\item If $X\not\cong\IP^1$ then $A$ is an embedding.
\item If $g\geq 2$ then $A(X)\subseteq J^{1,ss}_X$.
\item If $X$ is G-stable then $A(X)\subseteq J^{1,s}_X$.
\end{enumerate}
\end{prop}

\begin{prova}
The map $A$ is the composition of
$A_P\: X\to J^0_X$ followed by the
duality map $\lambda\:J^0_X\to J^0_X$, sending
$[\Ical]$ to $[\Ical^*]$, and the
translation $\tau\: J^0_X\to J^1_X$ by
$P$, sending $[\Ical]$ to $[\Ical\ox\o_X(P)]$. 
Assume $X\not\cong\IP^1$;
then
$A_P$ is an embedding by
Proposition~\ref{pointabel},
and since $\lambda$ and $\tau$ are
isomorphisms, also $A$ is an embedding.

Let us now see that $A(X)\subseteq J^{1,ss}_X$ if $g\geq 2$.
We claim that $\deg(\w|_Y)\geq 0$
for every subcurve $Y\subseteq X$. Indeed, since the
degree is additive, we need only
check the claim when $Y$ is irreducible. Then
$\deg(\w|_Y)\geq 0$ by \eqref{dwY1}, because
$g_Y\geq 0$ and, by hypothesis, $\delta_Y\geq 2$.
Thus, since $\w$ has degree $2g-2$, and
since $g\geq 1$ and $\delta_Y\geq 2$,
\begin{equation}\label{dwY}
\deg(\w|_Y)\leq 2(g-1)\leq\delta_Y(g-1).
\end{equation}
Now, for each $Q\in X$, there is a
natural inclusion $\o_X\to\mathfrak m_Q^*$, where $\mathfrak m_Q$ is the 
sheaf
of ideals of $Q$. Thus
\begin{equation}\label{dwQ}
\deg_Y(\mathfrak m_Q^*)\geq 0
\end{equation}
for every subcurve $Y\subseteq X$, and hence \eqref{dwY} implies
that $\mathfrak m_Q^*$ is semistable.

Finally, suppose that $X$
is G-stable. Let $Y$ be
a proper subcurve of $X$. Because
of \eqref{dwY} and \eqref{dwQ}, for the inequality
      $$\deg_Y(\mathfrak m_Q^*)\geq
         \frac{\deg_Y(\w)}{2g-2}-\frac{\delta_Y}{2}$$
to be an equality we would need that
$\deg_Y(\w)=\delta_Y(g-1)$. Since $\delta_Y\geq 2$,
we would need that $\delta_Y=2$ and
$\deg_{Y'}(\w)=0$. But this is not possible because $\w$ is ample.
Thus $A(X)\subseteq J^{1,s}_X$.
\end{prova}

\begin{rmk}\rm\setcounter{equation}{0}
There are special cases where
$J^{1,ss}_X=J^{1,s}_X$, for instance if
\begin{equation}\label{integer?}
\frac{\deg_Y(\w)}{2g-2}-\frac{\delta_Y}{2}
\end{equation}
is not an integer for any proper subcurve $Y\subset X$.
This will be the
case when $X$ is G-stable and
$g$ is odd. Indeed, suppose \eqref{integer?}
is an integer. Using \eqref{dwY1},
we have
$$
\frac{\deg_Y(\w)}{2g-2}-\frac{\delta_Y}{2}=
\frac{(2-g)\delta_Y-2\chi(\o_Y)}{2g-2}.
$$
Thus, if $g$ is odd, $\delta_Y$ must be even, and thus
$(2g-2)$ divides $\deg_Y(\w)$. However, as we saw in the
proof of Proposition \ref{abel}, this implies
that $\deg_Y(\w)=0$ or $\deg_{Y'}(\w)=0$, which is not possible
if $\w$ is ample.

If $X$ is a stable curve, in the sense of Deligne and
Mumford, it follows
from \cite{capoed} Prop. 3.15 that
$J^{1,ss}_X=J^{1,s}_X$, unless $X=Y_1\cup Y_2$, where
$Y_1$ and $Y_2$ are connected proper subcurves
of the same genus intersecting at
an odd number of points.
\end{rmk}

\begin{rmk}\label{dualrmk}\setcounter{equation}{0}\rm
Since $X$ is assumed Gorenstein, the dualizing map
$$
\lambda\:J_X\lra J_X;\  \  \  \  [\Ical]\mapsto[\Ical^*]
$$
is well-defined and takes $J_X^d$ isomorphically onto
$J_X^{-d}$, for every integer $d$.
Furthermore, given any vector bundle $F$ on $X$, using duality, we have
$$
\lambda(J_F^{ss})=J_{F^{\dagger}}^{ss},\quad\lambda(J_F^P)=J_{F^{\dagger}}^P
\quad\text{and}\quad
\lambda(J_F^s)=J_{F^{\dagger}}^s,
$$
where $F^\dagger:=F^*\ox\w$. In particular,
since $\mu(E_d^\dagger|_Y)=\mu(E_{-d}|_Y)$ for
each integer $d$ and each subcurve $Y\subseteq X$, it follows that
$$
\lambda(J_X^{d,ss})=J_X^{-d,ss},\quad\lambda(J_X^{d,P})=J_X^{-d,P}\quad\text{and}
\quad\lambda(J_X^{d,s})=J_X^{-d,s}.
$$
Thus, we could have defined $A_P$ as sending $Q$ to
$[\Ical_{Q/X}^*\ox\o_X(-P)]$, or $A$ as sending $Q$ to
$[\Ical_{Q/X}]$. Apart from the fact that the latter map
would have $J_X^{-1}$ as target instead of $J_X^1$, all the
conclusions would remain the same.

Essentially, the same observation applies to the twisted Abel maps
to be defined in \ref{twptembpre} and \ref{twembpre}.
\end{rmk}

\section{Twisted Abel maps of degree 0}

\begin{stp}\setcounter{equation}{0}\rm
(\emph{Spines and tails}).
A \emph{tail} of $X$ is a proper subcurve $Z\subset X$ with
$\delta_Z=1$. If $Z$ is a tail, so is $Z'$, and the unique point
$N$ of $Z\cap Z'$ is a separating node. In this case,
we say that $Z$ and $Z'$ are the
tails attached to $N$, and that $N$ generates
$Z$ and $Z'$. Notice that a tail is connected
(because $X$ is). A tail is called a
\emph{$P$-tail} if it does not contain $P$.
We denote by $\TX$ the set of all tails of $X$ and
by $\TXP$ the set of all $P$-tails.

A connected subcurve $Y$ of $X$ is called a
\emph{spine} if every point in
$Y\cap\ol{X-Y}$ is a separating node. In this case,
each connected component $Z$ of $\ol{X-Y}$ is a tail intersecting
$Y$ transversally at a single point on the nonsingular loci of $Y$ and $Z$.

Let $Y$ be a subcurve of $X$. If a singular point of $Y$ is a
separating node of $X$, then it is also a separating node of $Y$.
Conversely, if $Y$ is a spine then a separating node of $Y$ is a
separating node of $X$. As a consequence, if a subcurve $Z$ of $Y$
is a spine of $X$, then $Z$ is a spine of $Y$; conversely,
if $Z$ is a spine of $Y$ and $Y$ is a spine of $X$, then
$Z$ is a spine of $X$.

If $Y$ is a nonempty proper union of spines of $X$, then
any connected component of
$Y$ or $\ol{X-Y}$ is a spine. Two intersecting
spines with no common component intersect
transversally at a separating node of $X$.

A $q$-tuple $\mathfrak Z:=(Z_1,\dots,Z_q)$ of spines covering $X$,
each two with no component in common,  
is called a \emph{spine decomposition} of $X$. If
$Y$ is a spine of $X$, then $Y$ and the connected components
of $\ol{X-Y}$ form a spine decomposition of $X$.
\end{stp}

The following two lemmas will be much used.

\begin{lema}\label{ZZ}\setcounter{equation}{0}
Let $Z_1$ and $Z_2$ be tails of the curve $X$.
Then
$$
\text{either}\quad Z_1\cup Z_2=X\quad\text{or}\quad
Z_1\cap Z_2=\emptyset\quad\text{or}\quad
Z_1\subseteq Z_2\quad\text{or}\quad
Z_2\subsetneqq Z_1.
$$
\end{lema}

\begin{prova} This is \cite{capoed} Lemma 4.3.
\end{prova}

\begin{lema}\label{ctype}\setcounter{equation}{0} Let
$\mathfrak Z:=(Z_1,\dots,Z_q)$ be a spine decomposition of $X$.
Then there is an isomorphism
$$
u\:J_X\lra J_{Z_1}\times\cdots\times J_{Z_q}
$$
sending $[I]$ to $([I|_{Z_1}],\dots,[I|_{Z_q}])$.
Furthermore, for each integer $d$,
$$
u(J_X^d)=\bigcup_{d_1+\cdots+d_q=d}
J_{Z_1}^{d_1}\times\cdots\times J^{d_q}_{Z_q}.
$$
\end{lema}

\begin{prova} This is \cite{CoE} Prop. 3.2.
\end{prova}

\begin{stp}\setcounter{equation}{0}\label{twists}\rm
(\emph{Twisters on tails.})
By Lemma \ref{ctype},
or \cite{capoed} Lemma 4.4, for each tail $Z$ of the curve $X$
there is a unique, up to isomorphism, line bundle
on $X$ whose restrictions to $Z$ and $Z'$ are
$\o_Z(-N)$ and $\o_{Z'}(N)$, where $N$ is the separating
node generating $Z$. Denote this line bundle by $\o_X(Z)$.
We call it a \emph{twister}.

For simplification,
for each formal sum $\sum a_ZZ$ of tails $Z$ with integer
coefficients $a_Z$, set
$$
\o_X(\textstyle\sum a_ZZ):=\bigotimes\o_X(Z)^{\ox a_Z}.
$$

If $Z$ is a tail of $X$ attached to the node $N$, and
$f\:{\mathcal X}\to S$ is a one-parameter regular smoothing of $(X,N)$
(a flat, projective morphism of schemes such that
$S$ has dimension one, $X=f^{-1}(s)$ for a 
nonsingular 
$s\in S$, and ${\mathcal X}$ is smooth 
at $N$),
then $Z$ is a Cartier
divisor of ${\mathcal X}$, satisfying $\o_{\mathcal X}(Z)|_X\cong\o_X(Z)$ while
$\o_{\mathcal X}(Z)|_{f^{-1}(t)}=\o_{f^{-1}(t)}$ for each $t\in S\smallsetminus s$. So,
$\o_X(Z)$, though nontrivial, is the
limit of a family of trivial sheaves.
\end{stp}

\begin{stp}\label{twptembpre}\setcounter{equation}{0}\rm
(\emph{Degree-$0$ twisted Abel maps}). Let $Q\in X$.
If $Q$ is not a separating node, let $\Mcal_Q$
be the sheaf of ideals $\mathfrak m_Q$ of $Q$.
Notice that $\Mcal_Q$ is simple.
If $Q$ is a separating node, let $Z$ be the
$P$-tail generated by $Q$ (so that $P\not\in Z$)
and let $\Mcal_Q$ be the
unique line bundle on $X$ such that
$$
\Mcal_Q|_Z\cong\o_Z(-Q)\quad\text{and}\quad
\Mcal_Q|_{Z'}\cong\o_{Z'}.
$$
(That $\Mcal_Q$ exists and is unique follows from Lemma \ref{ctype}.)
Again, $\Mcal_Q$ is simple.

Define a map
\begin{equation}
\label{wtP}
A_P\:X\lra J^0_X; \  \   \  Q\mapsto [\Ical_Q]
\end{equation}
where
\begin{equation}
\label{IQ}
\Ical_Q:=\Mcal_Q\ox\o_X(P)\ox\o_X(-\sum_{Z\in \TXP\,;\,Z\ni Q} Z),
\end{equation}
the sum running over all tails $Z$ of $X$ containing
$Q$ but not $P$. Since $\Mcal_Q$ is simple, so is $\Ical_Q$, and hence $A_P$
is well-defined. We call $A_P$
the \emph{degree-$0$ (twisted) Abel map} of $X$ with base $P$.
If $X$ has no separating nodes then
(\ref{wtP}) coincides with (\ref{AP}).
We will see in Theorem \ref{twptemb} that, in any case,
$A_P$ is a
morphism of schemes.
\end{stp}

\begin{lema}\label{biglema}\setcounter{equation}{0}
Keep the notation of
\ref{twptembpre}. Let $W$ be a spine of $X$, and define
$$
B\: W\to J^0_W;\  \  \  Q\mapsto [\Ical_Q|_W].
$$
Then $B$ is a well-defined map and the
following three statements hold:
\begin{enumerate}
\item\label{bl1} If $P\in W$, then $B$ is the degree-$0$ Abel map of
$W$ with base $P$.
\item\label{bl2} If $P\in W'$, then $B$ is the
degree-$0$ Abel map of $W$ with base $N$, where $N$ is the
unique point of $W\cap W'$ on the same connected component
of $W'$ as $P$.
\item\label{bl3} In any case, the isomorphism class of
$\Ical_Q|_{W'}$ does not depend on $Q\in W$.
\end{enumerate}
\end{lema}

\begin{prova} We will use induction on $\delta_W$. Suppose first
that $\delta_W=1$, i.e., that $W$ is a tail. Let $N$ be the
separating node generating $W$.

Let $Q\in W$. Let $Z_1,\dots,Z_n$ be the $P$-tails of $X$
containing $Q$. It follows from Lemma~\ref{ZZ} that either
$Z_i\subset Z_j$ or $Z_j\subset Z_i$ for each distinct $i$ and
$j$. Thus, we may assume that
$$
Z_1\subset Z_2\subset\cdots\subset Z_{n-1}\subset Z_n.
$$
By definition, $\Ical_Q=\Mcal_Q\ox\Kcal_Q$,
where $\Kcal_Q:=\o_X(P)\ox\o_X(-Z_1-\cdots-Z_n)$.

Suppose first that $P\in W'$. Then $W\in\TXP$.
As $W\ni Q$, we have that
$W=Z_i$ for a certain $i$. The tails $Z_1,\dots,Z_{i-1}$
are also tails of $W$;
in fact, they are all the $N$-tails of $W$ containing $Q$.
And $Z_i,\dots,Z_n$ are all
the $P$-tails of $X$ containing $W$. So
\begin{align*}
\Kcal_Q|_W\cong&\o_W(N)\ox\o_W(-\sum_{Y\in\mathcal T_N(W)\,;\,Y\ni Q}Y)\\
\Kcal_Q|_{W'}\cong&\o_{W'}(P)\ox
\o_X(-\sum_{Y\in \TXP\,;\,Y\supseteq W} Y)|_{W'}.
\end{align*}
Notice that $\Kcal_Q|_{W'}$ does not depend on $Q$.

Since $Q\in W$ and $P\in W'$, we have that
$\Mcal_Q|_{W'}\cong\o_{W'}$, and hence $\Ical_Q|_{W'}$ does not depend 
on $Q$.
In addition, if $Q=N$ or $Q$ is not a separating node of $X$,
then $Q$ is not a separating node of $W$, and
$\Mcal_Q|_W$ is the sheaf of ideals of $Q$ in $W$. On the other hand, if
$Q$ is a separating node of $X$ different from $N$, then
$Q$ is a separating node of $W$; and if $Y$ is the $N$-tail of $W$
generated by $Q$, then $Y$ is the $P$-tail of $X$ generated by $Q$, and 
hence
$$
\Mcal_Q|_Y\cong\o_Y(-Q)\quad\text{and}\quad\Mcal_Q|_{\ol{W-Y}}\cong\o_{\ol{W-Y}}.
$$
In any case, it follows that
$[\Ical_Q|_W]$ is the image of $Q$ under the degree-0 Abel map of
$W$ with base $N$. This finishes the proof of the lemma in the
case where $P\in W'$

Now, suppose $P\in W$. We claim that
either $Z_i\cap W'=\emptyset$, or
$Z_i\supseteq W'$, for each integer $i$.
Indeed, $Z_i$ and $W'$ are both
$P$-tails. If $Z_i\cap W'\neq\emptyset$,
then either $Z_i\supseteq W'$ or
$Z_i\subseteq W'$ by Lemma~\ref{ZZ}.
However, if $Z_i\subseteq W'$, since $Q\in Z_i$ and
$P\not\in W'$, we have that $W'=Z_j$ for some
$j\geq i$. But, since $Q\in W$ as well, it follows that $Q=N$ and
$i=j=1$. But then $Z_i=W'$, proving the claim. Notice from our reasoning 
above
that $W'=Z_i$ for a certain $i$ if and only if $Q=N$.

If $Z_i\cap W'=\emptyset$, then $Z_i$ is tail of $W$. And if
$Z_i\supsetneqq W'$, then $\ol{Z_i-W'}$  is a tail of $W$.
On the other hand,
let $Y$ be a tail of $W$. If $N\not\in Y$, then
$Y$ is a tail of $X$ as well,
with $Y\cap W'=\emptyset$. And if $N\in Y$, then
$Y\cup W'$ is a tail of $X$. Thus
$$
\Kcal_Q|_W\cong\o_W(P)\ox\o_W(-\sum_{Y\in\mathcal T_P(W)\,:\,Y\ni Q}Y)\ox
\Lcal|_W,
$$
and $\Kcal_Q|_{W'}\cong\Lcal|_{W'}$,
where $\Lcal:=\o_X(-W')$ if $Q=N$ and $\Lcal:=\o_X$ otherwise.

Notice that $\Mcal_Q|_{W'}\cong\o_{W'}$ unless $Q=N$, in which case
$\Mcal_Q|_{W'}\cong\o_{W'}(-N)$. In any case,
$\Ical_Q|_{W'}$ is trivial, whence independent from $Q$. In addition,
$\Mcal_Q|_W$ is the sheaf of ideals $Q$ in $W$, if $Q$ is not a 
separating node
of $X$.  On the other hand, if $Q=N$ then
$\Mcal_Q|_W=\o_W$ and $\Lcal|_W\cong\o_W(-N)$.
And if $Q$ is a separating node of $X$ different from $N$, then
$Q$ is a separating node of $W$; and if $Y$ is the $P$-tail of $W$
generated by $Q$, then either $Y$ or $Y\cup W'$ is
the $P$-tail of $X$ generated by $Q$, depending on whether $N$ is on $Y$ 
or not,
and whence
$$
\Mcal_Q|_Y\cong\o_Y(-Q)\quad\text{and}\quad\Mcal_{\ol{W-Y}}\cong\o_{\ol{W-Y}}.
$$
In any case, it follows that $[\Ical_Q|_{W}]$ is the image of $Q$ under
the degree-0 Abel map of $W$ with base $P$, finishing the proof of the 
lemma
when $W$ is a tail.

Now, suppose $\delta_W>1$. Let $Z$ be a
connected components of $W'$, and $Y:=\ol{X-Z}$.
Then $Y$ is a tail. By induction, the map
$$
C\:Y\to J^0_Y;\  \  \  Q\mapsto [\Ical_Q|_Y]
$$
is well-defined, and is the degree-0 Abel map of $Y$ with base $P$
if $P\in Y$, and base $N'$ if $P\not\in Y$, where
$N'$ is the separating node of $X$ generating $Y$. Also
the isomorphism class of $\Ical_Q|_Z$ does not depend on
$Q\in Y$.

Note that $W$ is a spine of $Y$ and $\# W\cap\ol{Y-W}=\delta_W-1$.
So, by induction, $B$ is well-defined.
Furthermore, the isomorphism class of $\Ical_Q|_{\ol{Y-W}}$ does not
depend on $Q\in W$. Since neither does
the isomorphism class of $\Ical_Q|_Z$, and $\ol{Y-W}$ does
not intersect $Z$, we obtain \eqref{bl3}.

By induction, if $P\in W$, and hence $P\in Y$, then
$B$ is the degree-0 Abel map of $W$ with base $P$. If
$P\not\in W$, there are two cases to consider: If $P\not\in Y$,
that is, if $P\in Z$, then
$B$ is the degree-0 Abel map of $W$ with base $N'$; and
if $P\in Y$, then $P$ belongs to a connected component of
$\ol{X-W}$ other than $Z$, and thus $B$ is the degree-0 Abel map
of $W$ with base the point of intersection of this component
and $W$.
\end{prova}

\begin{defn}\setcounter{equation}{0}\label{sepline}{\rm
A smooth rational component $L$ of
the curve $X$ is called a \emph{separating line} if $L$
is a spine.}
\end{defn}

\begin{teo}\label{twptemb}\setcounter{equation}{0}
The degree-$0$ Abel map $A_P$
{\rm (defined in (\ref{wtP}))}
is a morphism of schemes. Furthermore:
\begin{enumerate}
\item
$A_P(X)\subseteq J_X^{0,P}$.
\item If $X$
contains no separating lines, then
$A_P$ is an embedding.
\end{enumerate}
\end{teo}

\begin{prova}
If $X$ has no separating nodes,
we
proved the statements in Proposition~\ref{pointabel}.
So, let $W$ be a tail of $X$.
We will assume,
by induction, that the theorem holds for curves with less
separating nodes than $X$.
Assume, without loss of generality, that $P\in W'$.
Let $N$ be the separating
node generating $W$. Notice that the separating nodes of $X$
are $N$ and those of $W$ and $W'$. Thus $W$ and $W'$
have fewer separating nodes than $X$.

By Lemma \ref{biglema}, under the identification
$J_X=J_W\times J_{W'}$ given by Lemma \ref{ctype},
\begin{equation}\label{A1B1}
A_P|_W=(A_1,B)\: W\lra J^0_W\times J^0_{W'},
\end{equation}
where $A_1$ is the degree-0 Abel map of $W$ with base $N$
and $B$ is constant, and
\begin{equation}\label{B2A2}
A_P|_{W'}=(B',A_2)\: W'\lra J^0_W\times J^0_{W'},
\end{equation}
where $A_2$ is the degree-0 Abel map of $W'$ with base $P$ and
$B'$ is constant. By induction, $A_P|_W$ and
$A_P|_{W'}$ are morphisms of schemes. Then so is $A_P$, because
$W$ and $W'$ intersect transversally at nonsingular points.

Now, let
$Q$ be a point of $X$. Let
$Z_1,Z_2\dots,Z_n$ be the $P$-tails of $X$ containing
$Q$. Using Lemma \ref{ZZ}, as in the proof of Lemma \ref{biglema},
we may assume that
$$
Z_1\subset Z_2\subset\cdots\subset Z_{n-1}\subset Z_n.
$$
Let
$\Kcal:=\o_X(-\sum Z_i)$.
Keep the notation of
\ref{twptembpre}. We want to show that $\Ical_Q$ is
$P$-quasistable. Let $Y$ be a connected proper subcurve of $X$.
It is enough to show that either $\deg_Y(\Ical_Q)\geq 0$, or
\begin{equation}\label{-1>-}
P\not\in Y\quad\text{and}\quad
\deg_Y(\Ical_Q)\geq -1\geq-\delta_Y/2.
\end{equation}

By \cite{capoed} Lemma 4.8, we have
$\deg_Y(\Kcal)\geq -1$.
Suppose first that
$\deg_Y(\Kcal)= -1$.
By the same lemma, $Y\subseteq Z'_1$.

Suppose $Q\in Y$. Then $Q\in Y\cap Y'$ and
$Q$ is the separating node generating $Z_1$.
So $\Ical_Q|_Y\cong\o_Y(P)\ox\Kcal|_Y$
if $P\in Y$ and
$\Ical_Q|_Y\cong\Kcal|_Y$
otherwise. If $P\in Y$,
then $\deg_Y(\Ical_Q)=0$. On the other hand, if $P\not\in Y$,
then $\delta_Y\geq 2$ and $\deg_Y(\Ical_Q)=-1$; so \eqref{-1>-} holds.

Now, suppose $Q\not\in Y$. Then
$\deg_Y(\Ical_Q)=0$ if $P\in Y$ and
$\deg_Y(\Ical_Q)=-1$ if $P\not\in Y$. We need only
show now that, if $Y$ is tail, then $P\in Y$.
Indeed, if $Y$ is a tail, then
$Y$ is generated by the same node that generates a
$Z_j$, for a certain $j$, again by \cite{capoed} Lemma 4.8.
Since $Y\subseteq Z'_1$, we have $Y=Z'_j$, and hence
$Y\ni P$.

Suppose now that
$\deg_Y(\Kcal)\geq 0$.
Then $\deg_Y(\Ical_Q)\geq -1$. If $Q\not\in Y$ or
$P\in Y$, then $\deg_Y(\Ical_Q)\geq 0$.
Now, suppose $Q\in Y$ and $P\not\in Y$. If
$\delta_Y\geq 2$, then \eqref{-1>-} holds. On
the other hand, if $Y$ is a tail, then $Y=Z_j$ for a
certain $j$, and hence
$\deg_Y(\Kcal)=1$.
In this case, $\deg_Y(\Ical_Q)=0$.

At any rate, $\Ical_Q$ is $P$-quasistable. Since this holds
for each $Q\in X$, the map $A_P$ factors through $J^{0,P}_X$.

Finally, assume $X$ contains no separating lines. First,
observe that a separating node of $W$ is a separating node
of $X$. So, given a smooth, rational component $L$ of $W$,
since
$$
L\cap L'\subseteq(L\cap\ol{W-L})\cup\{N\},
$$
not all points of $L\cap\ol{W-L}$ are separating nodes
of $W$. By induction, $A_1$ is an embedding.
By the same reasoning, $A_2$ is also an embedding.  
Hence $A_P|_{W}$ and $A_P|_{W'}$ are
embeddings.

Now, let $Q_1\in W$ and $Q_2\in W'$ and assume that
$A_P(Q_1)=A_P(Q_2)$. Then
$$A_1(Q_1)=B'(Q_2)\quad\text{and}\quad B(Q_1)=A_2(Q_2).$$
Since $A_1$ is injective, and $B'(W')=\{A_1(N)\}$,
we have that $Q_1=N$. Also, since $A_2$ is injective and
$B(W)=\{A_2(N)\}$, we have $Q_2=N$. Hence $Q_1=Q_2$. It follows
that $A_P$ is injective.

Also, since
$A_P|_W$ and $A_P|_{W'}$ are immersions, so
is $A_P$ everywhere but possibly at $N$.
But $A_P$ is an immersion also at $N$, because
$A_P|_W$ and $A_P|_{W'}$ are immersions at $N$,
and, under the identification $J_X=J_W\times J_{W'}$
given by Lemma \ref{ctype}, they take the tangent spaces of
$W$ and $W'$ at $N$ into
the linearly independent subspaces
$T_{J_W,A_1(N)}\oplus 0$ and
$0\oplus T_{J_{W'},A_2(N)}$ of $T_{J_X,A_P(N)}$,
respectively.

Thus $A_P$ separates points and tangent vectors.
Since $X$ is complete, and $A_P$ factors through
$J_X^{0,P}$, which is separated,
$A_P$ is an embedding.
\end{prova}

\section{Twisted Abel maps of degree 1}

\begin{stp}\label{smalltails}\setcounter{equation}{0}\rm
(\emph{Small tails and splitting nodes.})
We set now a rule that associates to every separating node $N$ of
the curve $X$ exactly one of the two tails that $N$ generates;
we shall call the chosen tail the \emph{small tail}
generated by $N$, and denote it by
$Z_N$. To do this, let $Z $ and $Z '$ be the two tails
generated by $N$; then $g_{Z}+g_{Z'}=g$.
There are two cases:
\begin{enumerate}
\item If $g_{Z}\neq g_{Z'}$, let $Z_N$ be the one
between $Z$ and $Z'$ having smaller genus. Thus $g_{Z_N}<g/2$.
\item\label{splitnode} If $g_{Z}= g_{Z'}=g/2$, make an arbitrary
choice between $Z$ and
$Z'$, and set it equal to
$Z_N$. In this case we call $N$ a {\it splitting node} of $X$.
\end{enumerate}

We denote by $\STX\subset\TX$ the set of all small tails of $X$.
By definition, there is
a  bijection between the set of separating nodes of $X$
and $\STX$. Observe that, if $X$ has a splitting node, $\STX$
depends upon the choice made in \eqref{splitnode} above.
\end{stp}

\begin{stp}\label{twembpre}\setcounter{equation}{0}\rm
(\emph{Degree-$1$ twisted Abel maps}).
For each  $Q\in X$
define a torsion-free, rank-1 sheaf $\Ncal_Q$ on $X$ as follows.
If $Q$ is not a separating node, let
$\Ncal_Q$
be the sheaf
of ideals of $Q$.
If $Q$ is a separating node,
let $Z\in \STX$ be the unique small tail attached to it, and
let $\Ncal_Q$ be the
unique line bundle on $X$ such that
$$
\Ncal_Q|_Z\cong\o_Z(-Q)\quad\text{and}\quad \Ncal_Q|_{Z'}\cong\o_{Z'}.
$$
(That $\Ncal_Q$ exists and is unique follows from Lemma \ref{ctype}.)
Note that $\Ncal_Q$ is simple.

Define a map
\begin{equation}\label{Atwist}
A\:X\lra J^1_X
\end{equation}
sending a point $Q$ of $X$ to $[\Ical^1_Q]$, where
\begin{equation}
\label{IQ1}
\Ical_Q^1:=(\Ncal_Q)^*\ox\o_X(\sum_{Z\in\STX\,:\, Z\ni Q} Z),
\end{equation}
the sum running over all small tails of $X$ containing $Q$.
Since $X$ is Gorenstein, and $\Ncal_Q$ is simple,
so is $\Ical_Q^1$, and hence
$A$
is well-defined. We call $A$
the \emph{degree-$1$ Abel map} of $X$. If $X$ has no separating
nodes then \eqref{Atwist} coincides with \eqref{Anotwist}. But
if $X$ admits a splitting node $N$,
then the definition of $A$ depends on the choice of the small tail
$Z_N$ associated to $N$ (see \ref{smalltails}). We will see in
Theorem \ref{twemb} that $A$ is a morphism of schemes.
\end{stp}

\begin{lema}
\label{ZZg}\setcounter{equation}{0}
If the curve $X$ is G-stable, then the following statements hold:
\begin{enumerate}
\item
\label{Z1} If $Z_1$ and $Z_2$ are tails of $X$ with
$Z_1\subsetneqq Z_2$, then $g_{Z_1}<g_{Z_2}$.
\item
\label{Z2} The curve $X$ has at most one splitting  node
{\rm (defined in \ref{smalltails})}.
\item
\label{Z3}
There is a point $Q$ on the nonsingular locus of $X$ such that
$\STX=\mathcal T_Q(X)$.
\end{enumerate}
\end{lema}

\begin{prova}
We prove Statement~\eqref{Z1}. As $X$ is G-stable,
$\w$ is an ample line bundle. So
$$
2g_{Z_1}-1=\deg_{Z_1}(\w)<\deg_{Z_2}(\w)=2g_{Z_2}-1,
$$
and hence $g_{Z_1}<g_{Z_2}$.

As for Statement~\eqref{Z2}, suppose $X$ has a splitting node $Q$,
and let $Z$ be a tail it generates.
Thus $g_Z=g_{Z'}=g/2$.
By contradiction, assume $X$ has another splitting node, $N$;
we may assume $N\in Z$.
Then $N$ generates two tails of $Z$, one of which is a tail of $X$.
So $Z$ contains properly a tail of $X$ of genus $g/2=g_Z$,
contradicting Statement~\eqref{Z1}.

Finally, let us prove the
last statement.
We may assume that $X$ has a tail, hence a small tail,
hence a maximal small tail, $Z$.
Let $Q$ be a point on the nonsingular locus of $X$ lying on the
irreducible component of $Z'$ containing $Z\cap Z'$. We claim
that $\STX=\mathcal T_Q(X)$. Indeed, let $W$ be a tail of $X$.
Suppose first that $Q\in W$. Then $W\cap Z\neq\emptyset$ and
$W\not\subseteq Z$. If
$W\cup Z=X$ then $W\supseteq Z'$, and since
$Z'$ is not small, neither is $W$.
If instead $W\cup Z\neq X$, then $Z\subsetneqq W$,
by Lemma~\ref{ZZ}, and thus $W$ is again not small, by
the maximality of $Z$. Conversely,
suppose $W$ is not small. Then $W'$ is.
But then, as we have just proved, $Q\not\in W'$,
and so $Q\in W$.
\end{prova}

\begin{teo}\label{twemb}\setcounter{equation}{0}
The degree-$1$ Abel map $A$ {\rm (defined in (\ref{Atwist}))}
is a morphism of schemes. Furthermore:
\begin{enumerate}
\item If $X$ has no separating lines, then $A$ is an
embedding.
\item\label{G-stableA} If $X$ is G-stable and $P$ does not lie on
any small tail of $X$, then $A(X)\subseteq J^{1,P}_X$.
\item\label{G-stableB}
If $X$ is G-stable then $A(X)\subseteq J^{1,ss}_X$.
\end{enumerate}
\end{teo}

\begin{prova}
By Lemma \ref{ZZg}, there is a point $Q$ on the nonsingular locus
of $X$ such that
$\STX=\mathcal T_Q(X)$, or equivalently, such that $Q$ does not lie
on any small tail of $X$. So Statement \eqref{G-stableB} follows
from \eqref{G-stableA}. Also, $A$ is a morphism
of schemes because, as in the proof of
Proposition~\ref{abel}, it is the composition of $A_Q$,
the degree-0 Abel map of $X$ with base $Q$,
with the duality map and the translation-by-$Q$
map. Furthermore, if $X$ contains no separating
lines, since the dualizing and translation maps
are isomorphisms, and since $A_Q$ is an embedding
by Theorem~\ref{twptemb}, also $A$ is an embedding.

Let us prove Statement \eqref{G-stableA}. Let $Q\in X$, and let
$Z_1,Z_2\dots,Z_n$ be the small tails of $X$ containing
$Q$. Using Lemma \ref{ZZ}, as in the proof of Lemma~\ref{biglema},
we may assume that
$$
Z_1\subset Z_2\subset\cdots\subset Z_{n-1}\subset Z_n.
$$
By hypothesis, $P\not\in Z_n$.

Keep the notation of \ref{twembpre}.
Let $W$ be a connected proper subcurve of $X$.
We need to show that
\begin{equation}\label{qst}
\deg_W(\Ical_Q^1)\geq
\frac{\deg_W(\w)}{2g-2}-\frac{\delta_W}{2},
\end{equation}
with equality only if
$P\not\in W$.

Let
$\Kcal:=\o_X(\sum Z_i)$.
By \cite{capoed} Lemma 4.8, we have
$\deg_W(\Kcal)\geq -1$.
Suppose first that
$\deg_W(\Kcal)=-1$.
Then $\deg_W(\Ical_Q^1)\geq -1$.
Also, by the same lemma, there is a unique integer $j$ such
that the separating node generating
$Z_j$ is contained in $W\cap W'$. In addition, $W\subseteq Z_j$,
and in particular, $W\subseteq Z_n$. So
$P\not\in W$.
Since $W\subseteq Z_n$,
and since $Z_n$ is a small tail,
$$
\deg(\w|_W)\leq\deg(\w|_{Z_n})\leq  g-1.
$$

So, if $\delta_W\geq 3$ then \eqref{qst} holds. On the other hand,
suppose $\delta_W=1$. Since $W\cap W'$ contains
the separating node of $Z_j$, and $W\subseteq Z_j$, we must
have $W=Z_j$. In this case, $\deg_W(\Ical_Q^1)=0$, and
hence \eqref{qst} holds as well.

Now, suppose $\delta_W=2$. We need to show
that $\deg_W(\Ical_Q^1)\geq 0$. Since $W$ contains the
separating node generating $Z_j$, and $W\subseteq Z_j$,
we have that $W=\ol{Z_j-Y}$ for a certain tail
$Y$ of $X$ properly contained in $Z_j$. Since $Z_j$ is a
small tail, so is $Y$. If $Y=Z_i$ for a certain $i<j$, then
$W\cap W'$ would also contain the separating node generating
$Z_i$, which is not possible. So, since $Y$ is a small tail,
$Q\not\in Y$. Thus $Q\in W$.  If $Q\not\in W\cap W'$ then
$\deg_W(\Ncal_Q^*)=1$, and hence $\deg_W(\Ical_Q^1)=0$. On the
other hand, if $Q\in W\cap W'$, since $Q\not\in Y$, it follows
that $Q$ is the separating node generating $Z_j$. Then $j=1$ and
$\deg_W(\Ncal_Q^*)=1$ as well, implying that $\deg_W(\Ical_Q^1)=0$.

The upshot is that \eqref{qst} holds and
$P\not\in W$ if $\deg_W(\Kcal)=-1$.
Now, suppose
$\deg_W(\Kcal)\geq 0$.
Then $\deg_W(\Ical_Q^1)\geq 0$. If $W$ is not a tail, then
\eqref{qst} holds, and the inequality is strict because
$W\neq X$, and hence $\deg_W(\w)<2g-2$. Suppose now
that $W$ is a tail. There are two cases to consider: $Q\in W$ and
$Q\not\in W$.

Suppose first that $Q\in W$.
If $W\not\subseteq Z'_1$,
then $\deg_W(\Ncal^*_Q)=1$, and hence
$\deg_W(\Ical_Q^1)\geq 1$. In this case, Inequality \eqref{qst}
holds and is strict. On the other hand,
if $W\subseteq Z'_1$, since $Q\in Z_1$,
we get that $Q$ is a separating node and $W=Z'_1$.
In this case,
$$
\deg_W(\Ical_Q^1)=\deg_W(\Kcal)=1,
$$
and thus Inequality \eqref{qst} holds as well and is strict.

Now, suppose $Q\not\in W$. Then $W'\ni Q$. If
$P\in W$,
then $W'$ is
a $P$-tail,
and hence a small tail.
Since $W'\ni Q$, we have that $W'=Z_i$ for a certain integer $i$,
and hence $\deg_W(\Kcal)=1$. Again, Inequality \eqref{qst}
holds and is strict. On the other hand, if
$P\not\in  W$,
then $W$ is a small
tail, and hence $\deg_W(\w)\leq g-1$. Since $\deg_W(\Ical_Q^1)\geq 0$,
Inequality \eqref{qst} holds.

At any rate, $\Ical_Q^1$ is
$P$-quasistable. Since this holds for every $Q\in X$, it follows
that $A(X)\subseteq J^{1,P}_X$.
\end{prova}

\section{Properties of the Abel maps}\label{img}

\begin{lema}
\label{g=0}\setcounter{equation}{0}
The curve $X$ has genus $g=0$ if and only if every irreducible
component of $X$ is a separating line.
\end{lema}

\begin{prova} By induction on
the number of irreducible components of $X$. If $X$ is
irreducible, $g=0$ implies that $X$ is smooth; so the lemma
holds trivially.
Suppose now that $X$ is reducible, and let $L$ be an irreducible
component of $X$ and $Z_1,\dots,Z_n$ the connected
components of $L'$. Since $X$ is Gorenstein and $L$ is irreducible,
$L$ is a separating line if and only if $g_L=0$ and
$\text{length}(L\cap Z_i)=1$ for each $i=1,\dots,n$ or,
equivalently,
$g_L=0$ and
\begin{equation}\label{LZ_i}
\text{length}(L\cap Z_1)+\cdots+\text{length}(L\cap Z_n)=n.
\end{equation}

Consider the
cohomology sequence associated to the natural
exact sequence
\begin{equation}\label{cohseq}
0\longrightarrow\o_X\longrightarrow\o_L\oplus
\o_{Z_1}\oplus\cdots\oplus\o_{Z_n}\longrightarrow
\o_{L\cap Z_1}\oplus\cdots\oplus\o_{L\cap Z_n}
\longrightarrow 0.
\end{equation}
If $h^1(X,\o_X)=g=0$ then $g_L=h^1(L,\o_L)=0$ and
\eqref{LZ_i} holds, whence $L$ is a separating line.

The converse is immediate, as    every irreducible component of $X$ is isomorphic to $\IP ^1$ and
every singularity is a separating node (see Example 9.8 in \cite{cner}).
\end{prova}

\begin{defn}\label{disctree}\setcounter{equation}{0}\rm
A \emph{separating tree of lines} of
the curve $X$ is a spine of (arithmetic) genus 0.
Equivalently, using Lemma~\ref{g=0},
a separating tree of lines of $X$ is a connected union of
separating lines of $X$.
\end{defn}

\begin{teo}\label{cptilde}\setcounter{equation}{0}
Let $A$ and $A_P$ be the Abel maps of the curve
$X$. Assume $g>0$. Set $B:=A$ or $B:=A_P$, and let $\wt X:=B(X)$.
Then the following statements hold:
\begin{enumerate}
\item
\label{emb}  Let $S$
be the union of all separating
lines of $X$.
Let $Y_1,\ldots,Y_n$ be the connected components
of $S'$.
Then $B|_{Y_i}$ is an
embedding for each $i=1,\dots,n$.
\item
\label{inj}  For any two distinct
points
$Q_1,Q_2\in X$,
we have that
$B(Q_1)=B(Q_2)$ if and only if
$Q_1$ and $Q_2$ belong to the same
separating tree of
lines.
\item
\label{R}  Let $L$ be a maximal separating tree of lines of $X$, and let
$R\in J_X$ for which $B(L)=\{R\}$.
Let $N_1,\dots,N_\delta$ be the points of $L\cap L'$.
Then $\wt X$ has an ordinary  $\delta$-fold singularity at $R$,
with linearly independent tangent lines equal to the images of the
differentials $d_{N_i}B:T_{X,N_i}\to T_{J_X,R}$.
\item\label{S4}
$\wt X$ is a curve of arithmetic genus $g$.
\end{enumerate}
\end{teo}

\begin{prova} Assume first that $B:=A_P$.
For each $Q\in X$, let $\Ical_Q$ be the simple, torsion-free, rank-1
sheaf on $X$ such that $B(Q)=[\Ical_Q]$.

We prove Statement \eqref{emb}. First, observe that
each $Y_i$ is a spine, and contains no separating lines.
Since $Y_i$ is a spine, using Lemma \ref{ctype},
it is enough to show that the map
$$
B_i\:Y_i\to J_{Y_i}; \  \   \  Q\mapsto [\Ical_Q|_{Y_i}]
$$
is an embedding. But, by Lemma \ref{biglema}, the map $B_i$
is a degree-0 Abel map. And,
since $Y_i$ contains no separating lines,
Theorem \ref{twptemb} implies that $B_i$ is an embedding.

Consider Statement \eqref{inj},
keeping the notation of Statement \eqref{emb}.
Suppose first that $Q_1$ and $Q_2$
belong to the same separating tree of lines.
Since the tree is connected, to show that
$B(Q_1)=B(Q_2)$ we may assume that
$Q_1$ and $Q_2$ lie on the smooth locus of $X$ and on
the same separating line, $L$. But then
$\Ical_{Q_1}|_L\cong\Ical_{Q_2}|_L$ and
$\Ical_{Q_1}|_{L'}\cong\Ical_{Q_2}|_{L'}$, as is easily seen.
Since $L$ is spine,
Lemma~\ref{ctype} yields $\Ical_{Q_1}\cong\Ical_{Q_2}$.

Suppose
now that $Q_1$ and $Q_2$ do not belong to the same tree of separating 
lines. We
must prove that $B(Q_1)\neq B(Q_2)$.
As we have just seen that $B$ is constant
along separating lines,
we may assume that $Q_1,Q_2\in Y_1\cup\dots\cup Y_n$.
Let $L_1,\dots,L_m$ be the connected components of $S$.
Since $X$ is connected, there are a positive integer $t$,
and integers $i_1,\dots,i_t\in\{1,\dots,n\}$ and
$j_1,\dots,j_{t-1}\in\{1,\dots,m\}$ such that $Q_1\in Y_{i_1}$ and
$Q_2\in Y_{i_t}$, while $Y_{i_\ell}\cap L_{j_\ell}\neq\emptyset$
and $L_{j_\ell}\cap Y_{i_{\ell+1}}\neq\emptyset$
for each $\ell=1,\dots,t-1$. Choose $t$ minimum; then
$Y_{i_\ell}\neq Y_{i_t}$ for every $\ell<t$. We will
show that $B(Q_1)\neq B(Q_2)$ by induction on $t$.

If $t=1$ then $Q_1$ and $Q_2$
belong to the same $Y_i$, and hence
$B(Q_1)\neq B(Q_2)$ by Statement \eqref{emb}.
Suppose now that $t\geq 2$. And suppose that $B(Q_1)=B(Q_2)$,
by contradiction. In particular, using Lemma \ref{ctype}, we have
$\Ical_{Q_1}|_{Y_{i_t}}\cong \Ical_{Q_2}|_{Y_{i_t}}$. But,
since $Y_{i_t}$ is a spine,
it follows from Lemma \ref{biglema} that
$\Ical_{Q_1}|_{Y_{i_t}}\cong \Ical_Q|_{Y_{i_t}}$
for each $Q$ on the same connected component of $Y'_{i_t}$ as
$Q_1$.By connectedness, there will be such $Q$ on
$L_{i_{t-1}}\cap Y_{i_t}$.
Since $\Ical_Q|_{Y_{i_t}}\cong \Ical_{Q_2}|_{Y_{i_t}}$, and since
$B_{i_t}$ is an embedding, as we saw
in the proof of Statement~\eqref{emb}, it follows
that $Q_2\in L_{i_{t-1}}$. Since $L_{i_{t-1}}$
is a separating tree of lines, $B(Q_2)=B(M)$ for any chosen
$M\in Y_{i_{t-1}}\cap L_{i_{t-1}}$. Since
$B(Q_1)=B(M)$, it follows by induction that $M=Q_1$.
But then $Q_1$ and $Q_2$ are on the same
separating tree of lines, namely $L_{i_{t-1}}$, reaching
a contradiction.

Finally, we prove Statements \eqref{R} and \eqref{S4}.
We proceed by induction
on the number of separating nodes of $X$. If zero,
then \eqref{R} is vacuous and \eqref{S4} follows
from \eqref{emb}.

Now, let $L$ be a maximal tree of separating lines of $X$, and
$R\in J_X$ such that $B(L)=\{R\}$. Notice that $L\neq X$
because $g>0$. Let $Z_1,\dots,Z_n$ be
the connected components of $L'$ and $N_1,\dots,N_\delta$ the
points in $L\cap L'$, Since $L$ is a spine with genus $0$,
we have that $n=\delta$ and
\begin{equation}\label{sumg}
g=g_{Z_1}+\cdots+g_{Z_n}.
\end{equation}
Also, up to reordering the $Z_i$,
we may assume that $Z_i$ is generated by $N_i$ for each $i=1,\dots,n$.
Notice that, for each $i$, no separating line of
$Z_i$ contains $N_i$, because otherwise its union with
$L$ would be a tree of separating lines of $X$ larger than $L$
itself. Also, notice that each $Z_i$ has less separating nodes
than $X$.

Since $(L,Z_1,\dots,Z_n)$ is a
spine decomposition of $X$, by Lemma \ref{ctype}
there is a natural isomorphism
$$
u=(u_L,u_1,\dots,u_n)\: J_X\lra
J_L\times J_{Z_1}\times\cdots\times J_{Z_n}.
$$
For each $i=1,\dots,n$, let $\wt Z_i:=B(Z_i)$.
By Lemma \ref{biglema}, as $Q$ moves on $Z_i$, the images
$u_L(B(Q))$ and $u_j(B(Q))$, for $j\neq i$, remain constant.
Thus $\wt Z_i\cong u_i(B(Z_i))$.
However, by the same Lemma \ref{biglema},
the composition $u_i\circ B|_{Z_i}$ is a degree-0 Abel map.
So, by induction, $\wt Z_i$ has genus $g_{Z_i}$.

Now, since $N_i$ is not contained in a separating line of $Z_i$,
it follows from Statement~\eqref{emb} that
$d_{N_i}B$ 
is injective on $T_{Z_i,N_i}$, and hence $\wt Z_i$
is nonsingular at $R$, with $T_{\wt Z_i,R}=d_{N_i}B(T_{Z_i,N_i})$.
Since $B$ contracts $L$, by Statement \eqref{inj}, it follows that
$T_{\wt Z_i,R}$ is the whole image of $d_{N_i}B$.

Using $u$ to make the identification
$$
T_{J_X,R}=T_{J_L,u_L(R)}\oplus T_{J_{Z_1},u_1(R)}\oplus
\cdots\oplus T_{J_{Z_n},u_n(R)},
$$
we may view $T_{\wt Z_i,R}$ as a subspace of
$$
0\oplus 0\oplus \cdots\oplus 0\oplus T_{J_{Z_i},u_i(R)}
\oplus 0\oplus\cdots \oplus 0
$$
for each $i=1,\dots,n$.
Thus the $T_{\wt Z_i,R}$ are linearly independent
subspaces of $T_{J_X,R}$. Since $R$ lies on the
nonsingular loci of all the $\wt Z_i$, and the $\wt Z_i$
cover $\wt X$, it follows that $\wt X$ has a
$n$-fold singularity at $R$,
with the $T_{\wt Z_i,R}$ for tangent lines,
finishing the proof of
Statement \eqref{R}.

As for Statement \eqref{S4}, since the $\wt Z_i$ intersect
transversally at $R$, we have exact sequences of the form:
$$
0\to\o_{\tilde Z_1}(-R)\to\o_{\tilde X}\to
\o_{\tilde Z_2\cup\cdots\cup\tilde Z_n}\to 0,
$$
$$
0\to\o_{\tilde Z_2}(-R)\to
\o_{\tilde Z_2\cup\cdots\cup\tilde Z_n}\to
\o_{\tilde Z_3\cup\cdots\cup\tilde Z_n}\to 0,
$$
$$
\vdots
$$
$$
0\to\o_{\tilde Z_{n-1}}(-R)\to
\o_{\tilde Z_{n-1}\cup\tilde Z_n}\to\o_{\tilde Z_n}\to 0.
$$
Computing Euler characteristics, and using that $\wt Z_i$ has genus
$g_{Z_i}$ for $i=1,\dots,n$,
$$
\chi(\o_{\wt X})=-g_{Z_1}-g_{Z_2}\cdots-g_{Z_{n-1}}+1-g_{Z_n}=1-g,
$$
where the last equality is
\eqref{sumg}. Thus $\wt X$ has genus $g$.

Finally, choosing a nonsingular point $Q$ of $X$ away from
all small tails of $X$, we have that $A\: X\to J_X^1$ is obtained from
$A_Q$ by taking duals and then translating by $\o_X(Q)$.
Since these operations are isomorphisms,
all of the statements proved above for $A_Q$ hold for
$A$.
\end{prova}

\section{Abel maps to Seshadri's compactified Jacobian}

Recall that for each point $Q$ of the curve $X$
we defined the simple, torsion-free, rank-$1$, degree-$1$ sheaf
$\Ical_Q^1$, and let $A(Q):=[\Ical_Q^1]$; see \ref{twembpre}.
When $X$ has a splitting node $N$, the definition of $\Ical_Q^1$
and that of $A$ depend on the choice of one of the two
tails generated by $N$ (the small tail $Z_N$).

If $X$ is $G$-stable, Theorem \ref{twemb} says
that $\Ical_Q^1$ is semistable
with respect to the canonical $1$-polarization
$E_1$ (cf. \ref{Ed}).  Its Jordan--H\"older filtrations
(cf. \ref{Sequiv}) are easy to describe.

\begin{lema}\label{W0lema}\setcounter{equation}{0}
Assume the curve $X$ is G-stable, and  let $Q\in X$. Then
$\Ical_Q^1$ is stable if and only if $X$ admits no splitting node.
If $X$ has a splitting node $N$, and $Z_N$ is the small tail
generated by $N$, then
$$
\emptyset\subsetneqq Z_N\subsetneqq X
$$
is the unique Jordan--H\"older filtration of $\Ical_Q^1$.
\end{lema}

\begin{prova}
As mentioned before the statement, $\Ical_Q^1$ is semistable
by Theorem \ref{twemb}.
So, given a proper subcurve $Y\subset X$, we have
\begin{equation}\label{qst=}
\deg_Y(\Ical_Q^1)\geq \frac{\deg_Y(\w)}{2g-2}-\frac{\delta_Y}{2}.
\end{equation}
We must show that
equality holds in \eqref{qst=} if and only if there is a splitting
node $N$ and $Y=Z_N$. For this, we may assume that $Y$
is connected.

Since $X$ is $G$-stable,
$$
0<\frac{\deg_Y(\omega)}{2g-2}<1.
$$
Also, as seen in the proof of Theorem \ref{twemb},
a consequence of \cite{capoed} Lemma 4.8, we have
$\deg_Y(\Ical^1_Q)\geq-1$. Hence, since $\delta_Y$ is a positive
integer, equality in \eqref{qst=} may hold only in two cases:
\begin{enumerate}
\item $\deg_Y(\Ical_Q^1)=-1$, while $\deg_Y(\w)=g-1$ and $\delta_Y=3$.
\item $\deg_Y(\Ical_Q^1)=0$, while $\deg_Y(\w)=g-1$ and $\delta_Y=1$.
\end{enumerate}

We claim that the first case is not possible.
Indeed, suppose by contradiction that it occurs.
Since
$\delta_Y=3$, we have $\deg_Y(\omega)=2g_Y+1$,
which implies that $g_Y=g/2-1$.
But since $\deg_Y(\Ical_Q^1)=-1$, it follows that
$$
\deg_Y(\o_X(\sum_{Z\in\STX\,;\, Z\ni Q}Z))=-1.
$$
So, by \cite{capoed} Lemma 4.8, we have that $Y\subseteq Z$ for a
certain small tail $Z$. Now,
$Y\neq Z$ because $\delta_Y=3$. Then,
since $\w$ is ample, $\deg_Z(\w)>\deg_Y(\w)$.
On the other hand, since $Z$ is a small tail,
$\deg_Z(\w)\leq g-1$, and hence $\deg_Z(\w)\leq \deg_Y(\w)$,
reaching a contradiction.

Now, suppose the second case occurs. Then $Y$ is
a tail of genus $g/2$. So $X$ has a splitting node $N$, and
we must show that $Y=Z_N$.
Suppose by contradiction that $Y\neq Z_N$ or, in other words,
that $Y'$ is a small tail.
There are two cases to consider:
$Q\in Y$ and $Q\in Y'$. If $Q\in Y'$, then
$Y'$ is the largest small
tail containing $Q$, and hence $\deg_Y(\Ical_Q^1)=1$, a
contradiction.

Suppose now that $Q\in Y$. Since $\deg_Y(\Ical_Q^1)\neq 1$,
there is a
small tail $Z$ of $X$ containing $Q$. For each such $Z$,
we have that $Z\not\supseteq Y$, because otherwise
the smallness of $Z$ would imply that $Y=Z$,
and $Y$ is not small.
Also, $Z\not\supseteq Y'$, because otherwise
the smallness of $Z$ would
imply that $Z=Y'$, and thus $Q\in Y\cap Y'$. But in this case,
$Y'$ would be the
unique small tail containing $Q$, implying that
$\deg_Y(\Ical_Q^1)=1$, a contradiction.
By Lemma~\ref{ZZ}, the only possibility left is that
$Z\subsetneqq Y$. But since this must hold for each small tail $Z$
containing $Q$, we would get
$\deg_Y(\Ical_Q^1)=1$ again, a contradiction.

Conversely, let $N$ be a splitting node, and suppose
$Y=Z_N$. We must show that $\deg_Y(\Ical_Q^1)=0$.
If $Q\in Y$, then $Y$
is the largest small tail containing $Q$, and hence
$\deg_Y(\Ical_Q^1)=0$. On the other hand, if
$Q\not\in Y$, then all small tails containing $Q$ are
strictly contained in $Y'$, and hence
$\deg_Y(\Ical_Q^1)=0$ as well.
\end{prova}

\begin{teo}\label{seshsep}\setcounter{equation}{0}
Assume that the curve $X$ is G-stable.
Let $A\:X\to J^1_X$ be the degree-$1$
Abel map {\rm (cf. \ref{twembpre})}
and $\Phi^1\: J^{1,ss}_X\to U_X(1)$ the S-map
{\rm (cf. \eqref{Sdmap})}.
Then the following two statements hold:
\begin{enumerate}
\item $\Phi^1\circ A$ is independent of the
choice of a small tail of genus $g/2$.
\item $\Phi^1$ restricts to a closed embedding on $A(X)$.
\end{enumerate}
\end{teo}

\begin{prova} Suppose that $X$ admits a
splitting node $N$ (unique by Lemma~\ref{ZZg}). Let $Z$
be a tail generated by $N$. Then
$A$ depends on whether we choose $Z$ or $Z'$ as small tail or,
in other words, whether we set $Z_N:=Z$ or $Z_N=Z'$.
Let $Q\in X$. Let $\Ical$ (resp. $\Ical'$)
be the torsion-free, rank-1 sheaf on $X$ such that
$A(Q)=[\Ical]$ (resp. $A(Q)=[\Ical']$), when $Z_N=Z$
(resp. $Z_N=Z'$). By Lemma \ref{W0lema}, we have
$\mathfrak S(\Ical)=\mathfrak S(\Ical')=\{Z,Z'\}$. Also,
$\Ical'|_Z\cong\Ical|_Z\ox\o_Z(N)$ and
$\Ical'|_{Z'}\cong\Ical|_{Z'}\ox\o_{Z'}(-N)$, and hence
$$
\text{Gr}(\Ical)=\Ical|_Z\oplus(\Ical|_{Z'}\ox\o_{Z'}(-N))\cong
(\Ical'|_Z\ox\o_Z(-N))\oplus \Ical'|_{Z'}=\text{Gr}(\Ical').
$$
So $\Ical$ and $\Ical'$ are S-equivalent, and thus
$\Phi^1(A(Q))$ is independent of the choice of $Z_N$.

As for the second statement, we need only show that
$\Phi^1|_{A(X)}$ separates points and tangent
vectors.
For each $Q\in X$, let $\Ical_Q^1$ be the torsion-free, rank-1
sheaf on $X$ such that $A(Q)=[\Ical_Q^1]$.
To show that $\Phi^1|_{A(X)}$ separates points, we need to show that
for each $Q,Q'\in X$,
\begin{equation}\label{IQQ'}
\mathfrak S(\Ical_Q^1)=\mathfrak S(\Ical_{Q'}^1)
\  \text{and}\  
\text{Gr}(\Ical_Q^1)\cong\text{Gr}(\Ical^1_{Q'})\quad
\text{if and only if}\quad
\Ical_Q^1\cong \Ical^1_{Q'}.
\end{equation}

Now, if $X$ contains no splitting node
then Lemma~\ref{W0lema} asserts
that $\Ical_Q^1$ and $\Ical^1_{Q'}$ are stable,
and thus \eqref{IQQ'} holds trivially.
On the other hand, if $X$ has a splitting node $N$,
then Lemma \ref{W0lema} implies that
$\mathfrak S(\Ical_Q^1)=\{Z,Z'\}$ and
$$
\text{Gr}(\Ical_Q^1)=\Ical_Q^1|_Z\oplus(\Ical_Q^1|_{Z'}\ox\o_{Z'}(-N)),
$$
where $Z:=Z_N$. An analogous description holds for $\Ical^1_{Q'}$. So
$\mathfrak S(\Ical_Q^1)=\mathfrak S(\Ical_{Q'}^1)$, and
$\text{Gr}(\Ical_Q^1)\cong\text{Gr}(\Ical^1_{Q'})$ if and only if
$\Ical_Q^1|_Z\cong \Ical^1_{Q'}|_Z$ and
$\Ical_Q^1|_{Z'}\cong \Ical^1_{Q'}|_{Z'}$, which,
by Lemma~\ref{ctype}, occurs if and only if
$\Ical_Q^1\cong \Ical^1_{Q'}$. So \eqref{IQQ'} holds.

Finally, for any $Q\in X$, Lemma \ref{W0lema} implies that
$\mathfrak S(\Ical^1_Q)$
is a spine decomposition, whether $X$ admits a splitting node or
not. So, given any nonzero $v\in T_{J_X,A(Q)}$,
it follows from \cite{CoE} Lemma~3.11 and Prop. 4.3
that there is a
subscheme $\Theta\subseteq U_X(a,d)$ such that
$d\Phi_{A(Q)}(v)\not\in T_{\Theta,\Phi(A(Q))}$. Thus
$\Phi|_{A(X)}$ separates
tangent vectors.
\end{prova}

\begin{rmk}\label{0bad}\setcounter{equation}{0}\rm
In general, it is not true that $\Phi^0\:J_X^{0,ss}\to U_X(0)$
restricts to a closed embedding on $A_P(X)$.
For a simple example, suppose $X$ is a nodal curve with four
irreducible components: two of them, $X_1$ and $X_2$, smooth and
rational, meeting at two points, $N_1$ and $N_2$; the third,
$X_3$, nonrational,
meeting only $X_1$ at a single point $N_3$, and the fourth,
$X_4$, also nonrational, meeting only $X_2$ at a single
point $N_4$. From its description, $X$ is stable. Also, it contains
no separating lines, so $A_P$ is an embedding by
Theorem \ref{twptemb}. Suppose
further that $g_{X_3}=g_{X_4}$. Then the only small tails of $X$
are $X_3$ and $X_4$.

Suppose $P\in X_1$. For each $Q\in X$, let
$\Ical_Q$ be the torsion-free, rank-1 sheaf on $X$ such that
$A_P(Q)=[\Ical_Q]$ (cf. \ref{twptembpre}). Let
$Q\in X_2-\{N_1,N_2,N_4\}$. Given a
proper subcurve $Y$ of $X$, we have that
$\deg_Y(\Ical_Q)=-\delta_Y/2$ if and only if $Y=X_2\cup X_4$. Thus
$\mathfrak S(\Ical_Q)=\{X_1\cup X_3,X_2\cup X_4\}$ and
$$
\text{Gr}(\Ical_Q)=\o_{X_2\cup X_4}(-Q)\oplus\o_{X_1\cup X_3}(P-N_1-N_2).
$$
Now, since $X_2$ is rational, Lemma \ref{ctype} yields  
that, as $Q$ moves on $X_2-\{N_1,N_2,N_4\}$,
the isomorphism class of $\o_{X_2\cup X_4}(-Q)$ does not change.
Therefore, the composition $\Phi^0\circ A_P$ contracts $X_2$.
\end{rmk}

\end{document}